\documentclass{gtart}   
\usepackage{amssymb,amsmath,amscd,epsf,eucal,latexsym} 

\input gtmonout
\volumenumber{2}
\volumeyear{1999}
\volumename{Proceedings of the Kirbyfest}
\pagenumbers{35}{86}
\papernumber{3}
\received{18 September 1998}\revised{13 April 1999}
\published{17 November 1999}

\def\xrightarrow#1{\buildrel{#1\,}\over\to}  

\newtheorem{theorem}{Theorem}[section]
\newtheorem{lemma}[theorem]{Lemma}
\newtheorem{cor}[theorem]{Corollary}
\newtheorem{prop}[theorem]{Proposition}
 
\theoremstyle{definition} 
\newtheorem{defn}[theorem]{Definition} 
\newtheorem*{rem}{Remark} 
\newtheorem*{rems}{Remarks} 

\newtheorem{example}{Example} 
 
\DeclareMathOperator{\intr}{int}

\DeclareMathOperator{\im}{im}
\DeclareMathOperator{\id}{id}
\DeclareMathOperator{\diff}{Diff}
\DeclareMathOperator{\homeo}{Homeo}

\renewcommand{\epsilon}{\varepsilon}

\newcommand{\zf}{\ZZ_{f}}
\newcommand{\cz}{\CC_{\ZZ}}

\newcommand{\Cz}{\mathbf{C}_{\ZZ}}

\newcommand{\cf}{\CC_{f}}

\newcommand{\dz}{\DD_{\ZZ}}

\newcommand{\df}{\DD_{f}}
\renewcommand{\dh}{\DD_{h}}
\newcommand{\bd}{\partial} 
 
\renewcommand{\o}{\circ} 
 
\newcommand{\ra}{\rightarrow} 

\newcommand{\hra}{\hookrightarrow} 
\newcommand{\ol}{\overline} 
\newcommand{\wh}{\widehat} 
\newcommand{\trb}{\bd_{\pitchfork}} 
\newcommand{\tb}{\bd_{\tau}}
\newcommand{\CO}{C^{0}} 
\newcommand{\CI}{C^{1}}
\newcommand{\CII}{C^{2}}
\newcommand{\Ci}{C^{\infty}}
\newcommand{\CC}{\mathcal{C}} 
\newcommand{\FF}{\mathcal{F}}

\newcommand{\LL}{\mathcal{L}} 
 
\newcommand{\MM}{\mathcal{M}}

\newcommand{\BB}{\mathcal{B}}

\newcommand{\UU}{\mathcal{U}} 
\renewcommand{\SS}{\mathcal{S}} 
 
\newcommand{\DD}{\mathcal{D}}
\newcommand{\ZZ}{\mathcal{Z}}

\newcommand{\R}{\mathbb{R}}
\newcommand{\Z}{\mathbb{Z}}

\newcommand{\Varphi}{\Phi}
\newcommand{\wt}{\widetilde} 
\newcommand{\0}{\emptyset} 
\newcommand{\hr}{H^{1}(M;\R)}

\newcommand{\sm}{\smallsetminus} 
\renewcommand{\ss}{\subset} 
\newcommand{\sseq}{\subseteq} 

\begin{document}

\title{Foliation Cones}                    
\authors{John Cantwell\\Lawrence Conlon}                  
\address{Department of Mathematics, St. Louis University, St. Louis, MO 
63103\\\smallskip\\
Department of Mathematics, Washington University, St.
Louis, MO 63130}
\asciiaddress{Department of Mathematics, St. Louis University\\ St. Louis, MO 
63103\\
Department of Mathematics, Washington University\\ St.
Louis, MO 63130}

\email{cantwelljc@slu.edu, lc@math.wustl.edu}

\begin{abstract} 
David Gabai showed that disk decomposable knot and link complements
carry taut foliations of depth one.  In an arbitrary sutured
$3$--manifold $M$, such foliations $\FF$, if they exist at all, are
determined up to isotopy by an associated ray $[\FF]$ issuing from the
origin in $H^{1} (M;\R)$ and meeting points of the integer lattice
$H^{1} (M;\Z)$.  Here we show that there is a finite family of
nonoverlapping, convex, polyhedral cones in $H^{1} (M;\R)$ such that
the rays meeting integer lattice points in the interiors of these
cones are exactly the rays $[\FF]$.  In the irreducible case, each of
these cones corresponds to a pseudo-Anosov flow and can be computed
by a Markov matrix associated to the flow.  Examples show that, in
disk decomposable cases, these are effectively computable.  Our result
extends to depth one a well known theorem of Thurston for fibered
3-manifolds. The depth one theory applies to higher depth as well.
\end{abstract}
\asciiabstract{David Gabai showed that disk decomposable knot and 
link complements carry taut foliations of depth one.  In an arbitrary
sutured 3-manifold M, such foliations F, if they exist at
all, are determined up to isotopy by an associated ray [F] issuing
from the origin in H^1(M;R) and meeting points of the integer
lattice H^1(M;Z).  Here we show that there is a finite family of
nonoverlapping, convex, polyhedral cones in H^1(M;R) such that
the rays meeting integer lattice points in the interiors of these
cones are exactly the rays [F].  In the irreducible case, each of
these cones corresponds to a pseudo-Anosov flow and can be computed
by a Markov matrix associated to the flow.  Examples show that, in
disk decomposable cases, these are effectively computable.  Our result
extends to depth one a well known theorem of Thurston for fibered
3-manifolds. The depth one theory applies to higher depth as well.}


\primaryclass{57R30}                
\secondaryclass{57M25, 58F15}              

\maketitlepage

\section{Introduction}

By theorems of Waldhausen~\cite{wald} and Thurston~\cite{th:norm}, the
classification of fibrations $\pi \co M\ra S^{1}$ which are transverse to
$\bd M$ is reduced to a finite problem for compact
3--manifolds. Indeed, the fibrations $\FF$ correspond one--one, up to
isotopy, to certain ``fibered'' rays $[\FF]\ss H^{1} (M;\R)$.  More
generally, the isotopy classes of $C^{2}$ foliations $\FF$ without
holonomy correspond one--one to ``foliated'' rays $[\FF]=\{t
[\omega]\}_{t\ge0} $, where $\omega $ is a closed, nonsingular 1--form
defining a foliation isotopic to $\FF$~\cite{LB,cc:LB}.  These rays
fill up the interiors of a finite family of convex, polyhedral cones
subtended by certain top dimensional faces of the unit ball $B\ss
H^{1} (M;\R)$ of the Thurston norm. The foliated rays that meet
nontrivial points of the integer lattice $H^{1} (M;\Z)$ are the
fibered rays.

In this paper, we consider sutured 3--manifolds $(M,\gamma )$ that
admit taut, transversely oriented depth one foliations $\FF$. We will
decompose $\bd M=\tb M\cup\trb M$ so that the compact leaves of $\FF$
are the components of the ``tangential boundary'' $\tb M$ and $\FF$ is
transverse to the ``transverse boundary'' $\trb M$.  In the language
of sutured manifolds~\cite{ga1}, $\trb M=\gamma $ and $\tb M=R (\gamma
)$. The depth one foliation fibers $M_{0}=M\sm\tb M$ over $S^{1}$ with
noncompact fibers.  Typically, these sutured manifolds will result
from cutting the complement $E (\kappa )$ of a $k$--component link
$\kappa $ along a Seifert surface $S$. In the resulting sutured
manifold $(M_S(\kappa ),\gamma )$, two copies of $S$ make up $\tb
M_{S} (\kappa )$ and $\trb M_S(\kappa )$ consists of $k$ annuli. In
more general examples, $\trb M$ has toral and/or annular components.

In~\cite{cc:isotopy}, we showed that such depth one foliations
correspond one--one, up to isotopy, to ``depth one foliated'' rays
$[\FF]\ss H^{1} (M;\R)$. Here we will show that there is a finite
family of convex, polyhedral cones in $H^{1} (M;\R)$, with disjoint
interiors, such that the rays through integer lattice points in the
interiors of the cones are exactly the depth one foliated
rays. Examples at the end of the paper will illustrate the fact that
this cone structure is often effectively computable.

In~\cite{cc:surg}, we exhibited families of depth one knot complements
$E (\kappa )$ in which the foliation cones could be described by a
norm on $H^{1} (M_S(\kappa );\R)$, but our examples will show that such
a description is generally impossible.  Instead, the dynamical
properties of flows transverse to the foliations will be exploited in
analogy with Fried's determination of the fibered faces of the
Thurston ball~\cite{fried}.  To show that the number of cones is
finite, we will use branched surfaces in the spirit of Oertel's
determination of the faces of the Thurston ball~\cite{Oer:hom}.

Again all rays in the interiors of these cones correspond to taut
foliations $\FF$ having holonomy only along $\tb M$, but those not
meeting the integer lattice $H^{1} (M;\Z)$ will have everywhere dense
noncompact leaves.  We conjecture that the isotopy class of such a
foliation is also uniquely determined by the foliated ray $[\FF]$.

The following theorem is meant to cover both the case of fibrations
and that of foliations of depth~one.  Accordingly, the term ``proper
foliated ray'' replaces the terms ``fibered ray'' and ``depth~one
foliated ray'' in the respective cases.  Here and throughout the
paper, $H^{1} (M)$ denotes de~Rham cohomology and explicit reference
to the coefficient ring $\R$ is omitted.

\begin{theorem}\label{cones}
Let $(M,\gamma )$ be a compact, connected, oriented, sutured
$3$--manifold.  If there are taut, transversely oriented foliations
$\FF$ of $M$ having holonomy (if at all) only on the
leaves in $\tb M$, then there are finitely many closed, convex,
polyhedral cones in $H^{1} (M)$, called foliation cones, having
disjoint interiors and such that the foliated rays $[\FF]$ are exactly
those lying in the interiors of these cones. The proper foliated rays
are exactly the foliated rays through points of the integer lattice
and determine the corresponding foliations up to isotopy.
\end{theorem}

In the fibered case, $\tb
M=\0$ and the theorem is due to Waldhausen and Thurston.

\begin{rem}
It will be necessary to allow the possibility that the entire vector
space $\hr$ is a foliation cone, this happening if and only if
$M=S\times I$ is the product of a compact surface $S$ and a compact
interval $I$ (Proposition~\ref{prod}). This is the one case in which
the vertex $0$ of the cone lies in its interior.  The foliated class
$0$ will correspond to the product foliation and $\{0\}$ will be a
(degenerate) proper foliated ray.
\end{rem}

 The proof of the theorem is reduced to the hyperbolic case where the
Handel--Miller theory of pseudo-Anosov endperiodic homeomorphisms
pertains. (This theory is unpublished, but cf \cite{fe:endp}). 
Determining the pseudo-Anosov monodromy for one
foliation $\FF$ gives rise to symbolic dynamics from which the
parameters for the foliation cone containing $[\FF]$ are easily read.
In the case of foliations arising from disk decompositions~\cite{ga0},
if the disks can be chosen in $M$ from the start, this procedure is
quite effective.  Indeed, the disks of the decomposition typically
split up in a natural way into the rectangles of a Markov partition
associated to the pseudo-Anosov monodromy. This partition determines
a finite set of ``minimal loops'' transverse to $\FF$ which span the
``tightest'' cone of transverse cycles in $H_{1} (M)$. The dual of
this cone is the maximal foliation cone containing $[\FF]$.

The authors thank Sergio Fenley for explaining to us many details of
the Handel--Miller theory (Section~\ref{psA}) and for a key step in
the proof of Theorem~\ref{also}. Research by the first author was
partially supported by N.S.F. Contract DMS--9201213 and that of the
second author by N.S.F. Contract DMS--9201723.

\section{Higher depth foliations}

Before proving Theorem~\ref{cones}, we indicate briefly its pertinence
to taut foliations of finite depth $k>1$ and smoothness class at least
$\CII$. To avoid technical problems, we assume that $\trb M=\0$. The
$\CII$ hypothesis guarantees that all junctures are compact, hence
that $(M,\FF)$ is homeomorphic to a $\Ci$--foliated
manifold~\cite[Main Theorem, page 4]{cc:smth2}.

\begin{rem}
The concept of a ``juncture'' is explained
in~\cite[Section~2]{cc:smth2} and will have important use in this
paper.
\end{rem}

Let $\SS\ss M$ be the compact lamination consisting of all leaves of
$\FF$ on which other leaves accumulate. While $\SS$ can have
infinitely many leaves, no real generality is lost by assuming it only
has finitely many leaves.  Indeed, it is possible to ``blow down''
finitely many foliated interval bundles in $\FF$ to produce a
foliation with $\SS$ finite-leaved.  The new foliation still has all
junctures compact, hence can be taken to be of class $\Ci$.

There are infinitely many ways to complete $\SS$ to a depth $k$
foliation. More precisely, consider any one of the components $U$ of
$M\sm\SS$, an open, connected, $\FF$--saturated set which is fibered
over $S^{1}$ by $\FF|U$.  The completion of $U$ relative to a
Riemannian metric on $M$ is a (generally noncompact) manifold $\wh{U}$
with boundary on which $\FF$ induces a depth one foliation $\wh{\FF}$.
As in~\cite[Theorem~1]{dip}, we write $$
\wh{U} = K\cup V_{1}\cup\dots\cup V_{n},
$$ where $K$ is a compact, connected, foliated, sutured manifold
(called the ``nucleus'') and $V_{i}\cong B_{i}\times I$ is a
noncompact, connected, foliated interval bundle (called an ``arm''),
$1\le i\le n$. Here, $\trb K$ has exactly $n$ components $A_{i}=\trb
V_{i}$, $1\le i\le n$. The assumption that $\trb M=\0$ implies that
these components are annuli. By choosing $K$ sufficiently large, one
guarantees that the foliation of each $V_{i}$ is the product
foliation.  This last assertion is due to compactness of the
junctures. Of course, $\wh{\FF}|K$ is of depth one.  Each depth one
foliation of $K$ which is trivial (that is, a product) at $\trb K$
determines a depth one foliation of $\wh{U}$, trivial in the arms. The
depth one foliations of $K$ that are trivial at $\trb K$ will be
called $\trb$--trivial.

In order to classify the $\trb$--trivial depth one foliations of $K$,
one first replaces $K$ with the manifold $K'$ obtained by gluing a
copy of $D^{2}\times I$ to each annular component of $\trb K$. The
$\trb$--trivial depth one foliations of $K$ correspond bijectively to
the depth one foliations of $K'$ with sole compact leaves the
components of $\bd K'$. Furthermore, there is a canonical splitting $$
H^{1}_{c} (K\sm\trb K)=H^{1} (K')\oplus V,
$$ where $V$ is spanned by the Poincar\'e duals of the components of
$\trb K$.  Thus, the foliation cones in $H^{1} (K')$ can be viewed as
``foliation cones'' in $H^{1}_{c} (K\sm\trb K)$. Obviously, these are
not full dimensional in the latter space, but they classify the
$\trb$--trivial depth one foliations and will be called the
$\trb$--trivial foliation cones of $K$.

 In order to classify all depth one foliations of $\wh{U}$, one must
allow an infinite exhaustion $$ K_{0}\ss K_{1}\ss\cdots\ss
K_{r}\ss\cdots\ss\wh{U}$$ by the potential nuclei. An inductive limit
process then leads to a finite family of ``foliation cones'' $C\times
\R^{N} \ss H^{1}_{c} (\wh{U})$, where $C$ ranges over the
$\trb$--trivial foliation cones of $K_{0}$ and $0\le N\le\infty$. We
omit the details.

Finally, under the assumption that $\SS$ has finitely many leaves,
this analysis only needs to be carried out for finitely many
open, saturated sets $U$.


\section{Reducing the sutured manifold}\label{reducing}

Let $M$ be a compact, connected, sutured 3--manifold.  A depth one
foliation determines a foliated class $[\omega ]\in H^{1} (M)$
which is represented by a \textit{foliated form} $\omega \in A^{1}
(M_{0})$.  This is to be a closed, nonsingular $1$--form which blows
up at $\tb M$ in such a way that the foliation $\FF_{0}$, defined by
$\omega $ on $M_{0}$, can be completed to a foliation $\FF$ of $M$,
integral to a $\CO$ plane field, by adjoining the components of $\tb
M$ as leaves. We will say that $\omega $ ``blows up nicely'' at $\tb
M$. The depth one condition implies that this form has period group of
rank one.  More generally, foliated forms of higher rank define
foliations $\FF$ tangent to $\tb M$ and such that each leaf of
$\FF|M_{0}$ is dense in $M$ and is without holonomy.  It can be shown
that all smooth foliations of $M$ having holonomy only along the
boundary leaves are $C^{0}$ isotopic to foliations defined by foliated
forms (Corollary~\ref{isot}). In case $M\cong S\times I$, we also
allow exact foliated forms.

In order to prove Theorem~\ref{cones}, it will be necessary to
``completely reduce'' $M$.  If $T\ss M$ is a compact, properly
imbedded surface, let $N(T)$ be a closed, normal neighborhood of $T$
in $M$ and denote by $N_{0}(T)$ the corresponding open, normal
neighborhood of $T$.  If $T_{1}\text{ and }T_{2}$ are disjoint,
properly imbedded surfaces, we always choose $N(T_{1})$ and $N(T_{2})$
to be disjoint.

\begin{defn}
Let $T\ss M$ be a properly imbedded, incompressible torus or annulus.
If $T$ is an annulus, require that one component of $\bd T$ lie on an
inwardly oriented component of $\tb M$ and the other component of $\bd
T$ on an outwardly oriented one.  Then $T$ is a reducing surface if it
is not isotopic through surfaces of the same type to a component of
$\trb M$.
\end{defn}

If $T\ss M$ is a reducing surface, we regard $M'=M\sm N_{0}(T)$ as a
(possibly disconnected) sutured manifold, $\trb M'$ being the union of
$\trb M$ and the two copies of $T$ in $\bd M'$.

\begin{defn}\label{compl:red}   If $T\ss M$ is a reducing surface, then 
$\{T\}$ is a reducing family.  Inductively, a reducing family is a
finite collection $\{T_{1},\dots,T_{r}\} $ of disjoint reducing
surfaces such that $\{T_{1},\dots,T_{r-1}\}$ is a reducing family and
$T_{r}$ is a reducing surface in whichever component of $M\sm
\bigcup_{i=1}^{r-1}N_{0}(T_{i})$ it lies. A maximal reducing family is
called a completely reducing family.  If $M$ contains no reducing
family, it is said to be completely reduced.
\end{defn}

If $M$ contains a reducing family, it contains a completely reducing
family.  In the following sections, we are going to prove
Theorem~\ref{cones} for the completely reduced case. Here, we will
show that this is sufficient.

\begin{theorem}\label{reduced}
If the conclusion of {\em Theorem~\ref{cones}} is true for
completely reduced, sutured $3$--manifolds, it is true for arbitrary
connected, sutured $3$--manifolds $M$.
\end{theorem}

Theorem~\ref{reduced} is proven by induction on the number $r$ of
elements of a completely reducing family in $M$. If $r=0$, there is
nothing to prove.  The inductive step is given by the following
lemmas.

\begin{lemma}\label{2 components} Let $T\ss M$ be a reducing surface,
$M'=M\sm N_{0}(T)$.  If $M'$ has two components and if the conclusion
of {\em Theorem~\ref{cones}} holds for each of these components, it
holds for $M$.
\end{lemma}

\begin{proof}
Indeed, fix an identification $N(T)=T\times [-1,1]$ and
 let $M_{+}'$ (respectively, $M_{-}'$) be the component of $M'$
meeting $N(T)$ along $T\times \{1\}$ (respectively, $T\times \{-1\}$).
Set $M_{\pm}=M'_{\pm}\cup N(T)$. Then \begin{align*}
M_{-}\cup M_{+} &= M\\ M_{-}\cap M_{+} &= N(T)
\end{align*} and Mayer--Vietoris gives an exact sequence
$$
0 \ra
\hr \xrightarrow{i} H^{1}(M_{-})\oplus H^{1}(M_{+}) \xrightarrow{j}
H^{1}(N(T)).
$$
Here, we use the
conventions that \begin{align*} i([\omega ]) &= ([\omega
|M_{-}],[\omega |M_{+}]) \\ j([\alpha ],[\beta ]) &= [\alpha
|N(T)]-[\beta |N(T)].  \end{align*} If $\omega $ is a foliated form,
the fact that $T$ is incompressible allows us to assume that $\omega
\pitchfork T$~\cite{rous:plonge,th:norm} so that $i ([\omega ])$
is a pair of foliated classes. By the inductive hypothesis, $i[\omega
]\in
\intr(\CC_{-}\times \CC_{+})$ for foliation cones $\CC_{\pm}$ in
$H^{1} (M_{\pm})$. Conversely, the Mayer--Vietoris sequence implies
that every class $[\omega ]$ carried into $\intr(\CC_{-}\times
\CC_{+})$ is represented by a foliated form obtained by piecing
together foliated forms $\omega _{\pm}$ representing classes in $\intr
(\CC_{\pm})$.  In order to make the cohomologous forms $\omega
_{\pm}|N (T)$ agree, one uses a theorem of Blank and
Laudenbach~\cite{LB}. It follows that the connected components of the
set of foliated classes in $\hr$ are exactly the interiors of a family
of convex, polyhedral cones of the form
$i^{-1}(\CC_{-}\times\CC_{+})$.
\end{proof}

We turn to the case that $M'$ is connected.  Again identify
$N(T)=N\times [-1,1]$, but realize $M'$ as $M\sm
\{T\times (-1/2,1/2)\}$. Set \begin{align*} N_{-} &= T\times
[-1,-1/2]\\ N_{+} &= T\times [1/2,1] \end{align*} and note that
\begin{align*} M'\cup N(T) &= M\\ M'\cap N(T) &= N_{-}\cup N_{+}.
\end{align*} The Mayer--Vietoris theorem gives the exact sequence
$$
\hr \xrightarrow{i} H^{1}(M')\oplus H^{1}(N(T)) \xrightarrow{j}
H^{1}(N_{-})\oplus H^{1}(N_{+}).
$$
One easily checks that the kernel of $i$ is 
spanned by the class that is Poincar\'e dual to $[T]\in H_{2}(M,\bd
M)$ and arguments analogous to the above prove the following.

\begin{lemma}\label{1 component}
If $M'$ is connected and the conclusion of {\em Theorem~\ref{cones}}
holds for $M'$, then it holds for $M$.
\end{lemma}

Theorem~\ref{reduced} follows. We turn to some further simplifying
conditions.

\begin{lemma}\label{ann-tor}
In {\em Theorem~\ref{cones}}, no generality is lost in assuming that
no component of $\tb M$ is an annulus or torus.
\end{lemma}

\begin{proof}
Indeed, let $\FF$ be defined by a closed, nonsingular 1--form $\omega $ which
blows up nicely at $\tb M$.  By the well understood structure of foliation
germs along toral and annular leaves, any such leaves in the boundary can be
perturbed inwardly by an arbitrarily small isotopy to become transverse to
$\FF$. Equivalently, $\omega $ is replaced by a nonsingular, cohomologous
form, differing from $\omega $ only in small neighborhoods of these boundary
leaves and transverse to them.  The former toral leaves are now components of
$\trb M$ and the former annular leaves are incorporated into transverse
boundary components. The foliated classes are unchanged.
\end{proof}

\begin{prop}\label{prod}
The sutured manifold $M$ has the form $S\times I$ if and only if every
class in $H^{1} (M)$ is a foliated class. That is, $H^{1} (M)$
is an entire foliation cone.
\end{prop}

\begin{proof}
Suppose that $M=S\times I$ and identify the compact interval as
$I=[-1,1]$.  Let $\lambda \co [-1,1]\ra[0,1]$ be smooth, strictly
positive on $(-1,1)$, and $\Ci$--tangent to $0$ at $\pm1$.  Each class
$[\gamma] \in H^{1} (S\times I)$ can be represented by a closed
form $\gamma =\omega +\lambda (t)^{-1} dt$ on $M_{0}=S\times (-1,1)$,
where $\omega $ is constant in the coordinate $t$ of $(-1,1)$.  This
form is nonsingular and blows up at the boundary.  The normalized form
$\lambda (t)\omega +dt$, while not closed, is defined and integrable
on all of $M$, determines the same foliation as $\gamma $ on $M_{0}$
and has the two components of $S\times \{\pm1\}$ as compact leaves.
If $\omega $ is not exact, this foliation has nontrivial holonomy
exactly on these boundary leaves, defining a foliation of the type we
are studying.  Note that $\gamma $ is exact if and only if the form
$\omega $ is exact, in which case the boundary leaves also have
trivial holonomy.  Thus, Reeb stability, coupled with Haefliger's
theorem~\cite{hae} that the union of compact leaves is compact,
implies that the foliation is isomorphic to the product foliation.  In
any event, the entire vector space $H^{1} (M)$ is the unique
foliation cone.  Conversely, suppose that the entire vector space is a
foliation cone. In particular, $0$ is a foliated class.  It is clear
that the foliations having nontrivial holonomy exactly on $\tb M$
correspond to nontrivial foliated classes, so Reeb stability and
Haefliger's theorem again imply that the foliation corresponding to
$0$ must be isomorphic to the product foliation on a manifold of the
form $S\times I$.
\end{proof}

\begin{prop}\label{tbnotempty}
In {\em Theorem~\ref{cones}}, it can be assumed without loss of
generality that $\tb M\ne\0$.
\end{prop}

Indeed, if $\tb M=\0$, the foliations we are studying are without
holonomy, the case of Theorem~\ref{cones} already covered by the
results of Waldhausen and Thurston.  In summary:

\begin{theorem}\label{ongoing}
The proof of {\em Theorem~\ref{cones}} is reduced to the case that $M$
is completely reduced, $\tb M\ne\0$ has no toral or annular components
and $M$ is not a product $S\times I$.
\end{theorem}

The hypotheses in Theorem~\ref{ongoing} will now be fixed as the
ongoing hypotheses in this paper.


\section{The  transverse structure cycles}\label{structurecycles}

Let $\LL$ be a 1--dimensional foliation of $M$, integral to a
nonsingular $\CO$ vector field (``leafwise $\CI$'') which is
transverse to $\tb M$ and tangent to $\trb M$. On the annular
components of $\trb M$, it is assumed that $\LL$ induces the product
foliation by compact intervals.  Let $\ZZ$ be the union of those
leaves of $\LL$ which do not meet $\tb M$.  It is evident that $\ZZ$
is a compact, 1--dimensional lamination of $ M_{0}$.  This lamination
is nonempty.  Indeed, if some leaf $\ell$ of $\LL$ issues from one
component of $\tb M$ but never reaches another, the asymptote of
$\ell$ in $M$ will be a nonempty subset of $\ZZ$.  Thus, if $\ZZ=\0$,
every leaf of $\LL$ issues from a component of $\tb M$ and ends at
another, implying that $M\cong S\times I$.  This contradicts one of
our ongoing hypotheses (Theorem~\ref{ongoing}). We call $\ZZ$ the
``core lamination'' of $\LL$. The following useful observation is left
as an exercise. 

\begin{lemma}\label{smooth_at_bd}
The foliation $\LL$ can be modified in a neighborhood of $\tb M$,
leaving $\ZZ$ unchanged, so that $\LL|M_{0}$ is integral to a
continuous vector field $v$ which is smooth near $\tb M$ and extends
smoothly to $\tb M$ so as to vanish identically there.
\end{lemma}

We will apply the Schwartzmann--Sullivan theory of asymptotic structure
cycles~\cite{sch_cycles,sull:cycles} to the core lamination $\ZZ$.
For this, the fact that the leaves of $\ZZ$ are integral to a vector
field which is at least continuous on $M$ will be essential. It is not
clear that we can significantly strengthen this regularity condition
for the endperiodic, pseudo-Anosov flows that will be needed in the
next section.

Let $\mu $ be a transverse, bounded, nontrivial, holonomy invariant
measure on $\ZZ$. Since the 1--dimensional leaves of the core
lamination $\ZZ$ have at most linear growth, such a measure exists by
a theorem of Plante~\cite{plante:meas}.  In standard fashion, 1--forms
$\omega $ on $M\text{ or }M_{0}$ can be integrated against $\mu $ in a
well defined way.  In local flow boxes, one integrates $\omega $ along
the plaques of $\ZZ$ and then integrates the resulting plaque function
against $\mu $.  Using a partition of unity, one assembles these local
integrals into a global one which is well defined because of the
holonomy invariance of $\mu $.  The resulting bounded linear
functional
$$\mu \co A^{1}(M)\ra\R$$ is a ``structure current'' of $\ZZ$ in the
sense of Sullivan~\cite{sull:cycles}.  The structure currents of $\ZZ$
form a closed, convex cone with compact base in the Montel space of
all 1--currents on $M$.

For a current $\mu $, defined as above by a transverse invariant
measure, it is easily seen that $\mu (df)=0$, for all smooth functions
$f$. That is, $\mu $ is a \textit{structure cycle} for $\ZZ$. The proof
uses the fundamental theorem of calculus (a.k.a. Stokes's theorem) on
plaques of $\ZZ$ and makes essential use of the fact that $\ZZ$ does
not meet $\tb M$.  Sullivan has proven (\textit{op.\!\!\!  cit.}) that the
transverse invariant measures are exactly the structure cycles and
form a closed, convex subcone $\Cz$ in the cone of structure currents
of $\ZZ$. The natural map of closed currents to homology classes
carries $\Cz$ onto a closed, convex cone $\cz\ss H_{1}(M)$ with
compact base. Each ray of the possibly infinite dimensional cone $\Cz$
is mapped one--one onto a ray of $\cz$, but $\cz$ is only finite
dimensional.

\begin{defn}
The dual cone $\dz$ to $\cz$ consists of all $[\omega ]\in \hr$
such that $\zeta (\omega )\ge0$, $\forall\,\zeta \in\Cz$.
\end{defn}

Let $\FF$ be a taut, transversely oriented foliation of $M$, smooth in
$M_{0}$, integral to a $\CO$ 2--plane field on $M$, having the
components of $\tb M$ as sole compact leaves and having nontrivial
holonomy exactly on these compact leaves.  We say that the tautly
foliated, sutured manifold $(M,\FF)$ is ``almost without
holonomy''. Remark that foliations defined by foliated forms are of
this type, but the converse is not quite true.  By a theorem of
Sacksteder~\cite{sa:pseudo}, $\FF|M_{0}$ admits a transverse,
continuous, holonomy invariant measure $\nu $ which is finite on
compact sets. By our hypothesis that the leaves in $\tb M$ are exactly
the ones with nontrivial holonomy, such a nontrivial holonomy
transformation must have no fixed points in $M_{0}$, hence $\nu $
becomes unbounded near $\tb M$. One should think of $\nu $ as a
``$C^{0}$ foliated form''.

Suppose that the one dimensional foliation $\LL$ is transverse to
$\FF$. As is well known, $\nu $ has well
defined line integrals in $M_{0}$ and $\int_{\sigma }^{}d \nu $
depends only on the homology class of a loop $\sigma $.  Thus, $\nu $
is a cocycle and can play a role analogous to that of a closed
1--form, allowing us to integrate $\nu $ against the structure cycle
$\mu $ of $\ZZ$ in a well defined way.  This integral is clearly
positive, proving the ``only if'' part of the  following.

\begin{theorem}\label{dual cone} The one dimensional foliation $\LL$
is transverse to a foliation $\FF$ which is almost without holonomy if
and only if no nontrivial structure cycle of $\ZZ$ bounds. In this
case, the dual cone $\dz$ has nonempty interior and every element of
$\intr \dz$ is a foliated class, represented by a foliated form which
is transverse to $\LL$.
\end{theorem}

\begin{proof} We prove the ``if'' part and the subsequent assertions. 
By~\cite[Theorem I.7, part iv]{sull:cycles}, the hypothesis that no
nontrivial structure cycle bounds implies that $\intr \dz\ne \0$ and
that each class in this open cone has a representative form $\omega\in
A^{1}(M) $ which is transverse to $\ZZ$.  We must replace $\omega $
with a form $\wt{\omega }=\omega +dg$ such that $g\in
C^{\infty}(M_{0})$, $dg$ blows up nicely at $\tb M$ and $\wt{\omega
}\in A^{1}(M_{0})$ is transverse to $\LL|M_{0}$. Remark that $dg$
blows up nicely at $\tb M$ if $dg$ is unbounded in a neighborhood $V$
of $\tb M$ and the smooth foliation of $V\cap M_{0}$ by the level sets
of $g$ can be completed to a foliation of $V$ integral to a $\CO$
plane field by adjoining the components of $\tb M$ as leaves.  Since
$\omega $ is bounded on the compact manifold $M$, it is clear that
$\wt{\omega }$ will blow up nicely at the boundary, hence be a
foliated class, and that $[\omega ]=[\wt{\omega }]$.

Let $v$ be a $\CO$ vector field to which $\LL|M_{0}$ is integral and
which extends by $0$ to a continuous field on $M$.  By
Lemma~\ref{smooth_at_bd}, we assume that $v$ is smooth in a
neighborhood $V$ of $\tb M$.  Let $S_{-}$ be the union of the inwardly
oriented components of $\tb M$ and $S_{+}$ the union of the outwardly
oriented ones.  We can assume that the part of $V$ bordering $S_{-}$
is parametrized by the local $v$--flow as $S_{-}\times [-\infty,0)$
and the part bordering $S_{+}$ as $S_{+}\times (0,\infty]$. In
particular, the fibers $\{x\}\times [-\infty,0)$ (respectively,
$\{x\}\times (0,\infty]$) are subarcs of leaves of $\LL$,
$\forall\,x\in S_{-}$ (respectively, $\forall\,x\in S_{+}$). Let $W\ss
V$ be a neighborhood of $\tb M$ such that $\ol{W}\ss V$.

Define $f\in C^{\infty}(S_{-}\times (-\infty,0))$ so that
$f(x,t)=f(t)$ everywhere,  $$f(t) = 
\begin{cases}
t, & t\le-1,\\ 0, & -\frac{1}{2}\le t,
\end{cases}
$$ and so that $f'(t)\ge0$ on $(-\infty,0)$.  Define $f$ analogously
on $S_{+}\times (0,\infty)$ and extend these definitions by 0 to a
smooth function $f$ on all of $M_{0}$. Remark that $df(v)\equiv1$ in
$V\cap M_{0}$, hence the form $df$ clearly blows up nicely at $\tb M$.

There is an open neighborhood $U$ of $\ZZ$ in $M_{0}$ such that
$\omega (v)>0$ on $U$. Let $z\in M_{0}\sm\ZZ$ and let $\ell_{z}$
denote a compact subarc of the leaf of $\LL$ through $z$ which has $z$
in its interior and exactly one end in $\tb M$. For definiteness,
assume this end lies in $S_{-}$.  There are open neighborhoods
$V_{z}\text{ and }U_{z}$ of $\ell_{z}$, $\ol{V}_{z}\ss U_{z}$, and a
smooth field $w_{z}$ on $M$ which approximates $v$ arbitrarily well on
$V_{z}$, approximates the \textit{direction} of $v$ arbitrarily well on
$U_{z}\cup V$, agrees with $v$ on $W$ and vanishes identically on
$M\sm (U_{z}\cup V)$.  Let $\Varphi ^{z}_{t}$ be the flow of $w_{z}$.
If $x\in W\cap V_{z}$ near $S_{-}$ and $\tau_{z} \in\R$ are such that
the $z=\Varphi ^{z}_{\tau _{z}} (x)$, choose $t_{z}>\tau _{z}$ and let
$$ df_{z}=\Varphi _{-t_{z}}^{*} (df).  $$ Thus, we can assume that
$df_{z} (v)\ge0$ on $M_{0}$ and $>1/2$ on $(V_{z}\cup W)\cap M_{0}$.
By compactness of $M$ and boundedness of $\omega $, find finitely many
points $z_{i}\in M_{0}\sm\ZZ$, $1\le i\le r$, and a constant $c>0$
such that $\{U,V_{z_{i}}\}_{i=1}^{r}$ is an open cover of $M$ and $$
\wt{\omega }=
\omega +dg = \omega +c\sum_{i=1}^{r}df_{z_{i}}
$$ is a closed form, strictly positive on $v|M_{0}$ which blows up at
$\tb M$. In fact, we have guaranteed that $dg$ is a constant multiple
of $df$ near $\tb M$, so $\wt{\omega }$ blows up nicely at $\tb M$.
\end{proof}
 
Examples show that, if $\FF|M_{0}$ is dense-leaved, the invariant
measure $\nu $ may not be absolutely continuous, let alone smooth, in
which case $\FF$ is not defined by a foliated form. Although these
foliations are not the primary focus of this paper, we note the
following.

\begin{cor}\label{isot}
The foliation $\FF$ is $\CO$ isotopic to a foliation $\wt{\FF}$
defined by a foliated form $\wt{\omega} $.
\end{cor}

Indeed, $[\nu ]\in \dz$, so Theorem~\ref{dual cone} provides a
foliated form $\wt{\omega} $ representing this same class.  We take
$\wt{\FF}$ to be the foliation defined by $\wt{\omega} $ and find an
isotopy of $\FF$ along $\LL$ to $\wt{\FF}$, using the continuous
measures on the leaves of $\LL$ induced by $\nu $ and $\wt{\omega} $.
For details, see~\cite[Section~2]{cc:LB}.

It will be important to characterize a particularly simple spanning
set of $\Cz$, the so called ``homology directions'' of Fried~\cite[page
260]{fried}.  Assuming that $\LL|M_{0}$ has been parametrized as a
nonsingular $\CO$ flow $\Varphi _{t}$, select a point $x\in\ZZ$ and
let $\Gamma $ denote the $\Varphi $--orbit of $x$.  If this is a
closed orbit, it defines a structure cycle which we will denote by
$\ol{\Gamma }$.  If it is not a closed orbit, let $\Gamma _{\tau
}=\{\Varphi _{t} (x)\mid 0\le t\le \tau \}$.  Let $\tau
_{k}\uparrow\infty$ and set $\Gamma _{k}=\Gamma _{\tau _{k}}$.  After
passing to a subsequence, we obtain a structure current $$ \ol{\Gamma}
= \lim_{k\ra\infty}\frac{1}{\tau _{k}}\int_{\Gamma _{k}}^{}. $$

\begin{lemma}\label{long_aco}
A structure current $\ol{\Gamma }$, obtained as above, is a structure
cycle.
\end{lemma}

\begin{proof}
By compactness of $\ZZ$, we can again pass to a subsequence so as to
assume that the points $\Varphi _{\tau_{k} } (x)$ all lie in the same
flowbox $B$, $k\ge1$.  Thus, we can close up $\Gamma _{k}$ to a loop
$\Gamma _{k}^{*}$ by adjoining an arc in $B$ from $x_{k}$ to
$x$. These arcs can be kept uniformly bounded in length, hence the
sequence of \textit{singular cycles} $(1/\tau _{k})\Gamma _{k}^{*}$
(generally not foliation cycles) also converges to $\ol{\Gamma}
$. These approximating singular cycles are called ``long, almost
closed orbits'' of $\Varphi _{t}$. Since the space of cycles is a
closed subspace of the space of currents, this proves that
$\ol{\Gamma} $ is a cycle.
\end{proof}

\begin{defn}
All  structure cycles $\ol{\Gamma} $, obtained as
above, and  their homology classes are called homology directions of
$\ZZ$.
\end{defn}

An elementary application of ergodic theory proves the following
({cf}~\cite[Proposition II.25]{sull:cycles}).

\begin{lemma}\label{span}
Any structure cycle $\mu \in\Cz$ can be arbitrarily well approximated by
finite linear combinations $\sum_{i=1}^{r}a_{i}\ol{\Gamma }_{i}$ of homology
directions. If $\mu \ne0$, the coefficients $a_{i}$ are strictly positive and
their sum is bounded below by a constant $b_{\mu }>0$ depending only on $\mu
$.
\end{lemma}

While much of the Schwartzmann--Sullivan theory requires that $\LL$ be
at least leafwise $\CI$, this lemma suggests how to define the cones
$\cz\text{ and }\dz$, even for the case in which the transverse
foliation $\LL$ is only $\CO$. Indeed, the long, almost closed orbits
(and the honest closed orbits) are defined, they are singular homology
cycles, their classes in $H_{1} (M)$ form a bounded set and limit
classes as above, taken over sequences for which $\tau _{k}\uparrow
\infty$, are called homology directions.  The closure of the set of
positive linear combinations of homology directions forms a convex
cone $\cz\ss H_{1} (M)$ and $\dz\ss H^{1} (M)$ is the dual cone. While
results such as Theorem~\ref{dual cone} become problematic in this
context, we will always be in a position to identify these cones with
ones coming from leafwise $\CI$ data.

 In what follows, the foliation $\FF$ will be of depth one and the
transverse foliations $\LL\text{ and }\LL'$ may only be $\CO$.

\begin{lemma}\label{C0} Let $\LL\text{ and }\LL'$ be two $1$--dimensional
foliations transverse to $\FF$.  Suppose that the respective core
laminations $\ZZ\text{ and }\ZZ'$ are $\CO$--isotopic by an isotopy
$\varphi _{t}\co \ZZ\hra M$, $\varphi _{0}|
\ZZ=\id_{\ZZ}$ and $\varphi _{1} (\ZZ)=\ZZ'$, such that $\varphi _{t} (x)$
lies in the same leaf of $\FF$ for $0\le t\le 1$,
$\forall\,x\in\ZZ$. Then $\cz=\CC_{\ZZ'}$.
\end{lemma}

\begin{proof}  Parametrize the two foliations as  flows using the
same transverse invariant measure for $\FF$.  Since $\FF$ is
leafwise invariant under the  isotopy, the flow parameter is
preserved and the long, almost closed orbits of $\ZZ$ are isotoped to
the long, almost closed orbits of $\ZZ'$.  Homotopic singular cycles
are homologous and the assertions follow.
\end{proof}

Note that we do not require this to be an ambient isotopy. The
property that points of $\ZZ$ remain in the same leaf of $\FF$
throughout the isotopy will be indicated, as above, by saying that
$\FF$ is leafwise invariant by $\varphi _{t}$.

The following is proven using the well understood structure of depth
one foliations in neighborhoods of $\tb M$.

\begin{lemma}\label{1st:ret'}
If $\LL\text{ and }\LL'$ are two transverse foliations which induce
the same first return map $f$ on a noncompact leaf $L_{0}$ of $\FF$,
then, without changing this property, one can modify $\LL'$ in an
arbitrarily small neighborhood of $\tb M$ to agree with $\LL$ in a
smaller neighborhood.
\end{lemma}

\begin{lemma}\label{1st:ret}
Let $\LL\text{ and }\LL'$ be two $1$--dimensional foliations
transverse to $\FF$ and inducing the same first return map $f\co L_{0}\ra
L_{0}$ on some leaf $L_{0}$ of $\FF|M_{0}$.  Then $\ZZ'$ is isotopic
to $\ZZ$ by a $\CO$ isotopy leaving $\FF$ leafwise invariant, hence
$\CC_{\ZZ}=\CC_{\ZZ'}$ and $\dz=\DD_{\ZZ'}$.
\end{lemma}

\begin{proof}  Parametrize the two
foliations, using the same transverse invariant measure for $\FF$, so as to
obtain flows $\Varphi _{t}$ and $\Varphi '_{t}$, both carrying leaves $L$ of
$\FF$ to leaves $\Varphi _{t} (L)=\Varphi '_{t} (L)$ of $\FF$, and such that
$$\Varphi _{1}|L_{0}=f=\Varphi '_{1}|L_{0}.$$ We choose a compact, connected
submanifold $K\ss L_{0}$, separating all the ends of $L_{0}$, with $Z\ss\intr
K$.  By Lemma~\ref{1st:ret'}, we lose no generality in assuming that $\LL$ and
$\LL'$ coincide near $\tb M$, so we can choose $K$ larger, if necessary, to
guarantee that $\Varphi _{t} (\bd K)=\Varphi '_{t} (\bd K)$, $0\le t\le 1$.
These properties imply that $\Varphi _{t} (K)=\Varphi '_{t} (K)$, $0\le
t\le1$.  Define
\begin{gather*}
\varphi _{t}\co K\ra K,\quad 0\le t\le1,\\ \varphi _{t}=\Varphi '_{-t}\o
\Varphi _{t}|K,
\end{gather*}  a loop in $\homeo_{\,0} (K)$ based at  $\varphi
_{0}=\varphi _{1}=\id_{K}$. Here, $\homeo_{\,0} (K)$ denotes the
identity component of the group of homeomorphisms of $K$ (with the
compact--open topology). Since no end of $L_{0}$ has a neighborhood of
the form $\R\times I$ or $\R\times S^{1}$ (a consequence of our
ongoing hypotheses), we can choose $K$ to have negative Euler
characteristic. Theorems of M. E.
Hamstrom~\cite{ham:disk-holes,ham:torus,ham:homeo} then imply that the
group $\homeo_{\,0} (K)$ is simply connected (see the remark below), so
there is a homotopy $\varphi _{t}^{s}$ in $\homeo_{\,0} (K)$, $0\le
s\le1$, fixing the basepoint, with $\varphi ^{0}_{t}=\varphi _{t}$ and
$\varphi _{t}^{1}=\id_{K}$, $0\le t\le1$.  This defines a continuous
deformation of $$\Varphi _{t}|K=\Varphi '_{t}\o \varphi _{t}^{0}\text{
to }
\Varphi '_{t}\o
\varphi _{t}^{1}=\Varphi '_{t}|K$$ which slides points along the
leaves of $\FF$.  This restricts to a $\CO$ isotopy of $$
\ZZ\ss\bigcup_{0\le t\le1}\Varphi _{t} (K)\text{ to } \ZZ'\ss 
\bigcup_{0\le t\le1}\Varphi '_{t} (K) 
$$ and everything now follows by Lemma~\ref{C0}.
\end{proof}

\begin{rem}
The theorem of Hamstrom, cited in the above proof, is that the group
$\homeo_{\,0} (K,\bd K)$ which fixes $\bd K\ne\0$ pointwise is
homotopically trivial. This is true whether or not $\chi (K)$ is
negative.  The assumption of negative Euler characteristic implies
that $\pi _{1} (K,x)$ is free on at least two generators, in which
case one shows that any loop $\varphi $ on $\homeo_{\,0} (K)$, based at
the identity, is base point homotopic to a loop in the subgroup
$\homeo_{\,0} (K,\bd K)$.  It follows that $\homeo_{\,0} (K)$ is simply
connected.  Indeed, for each point $x\in K$, $\varphi $ defines a loop
$\varphi _{x}$ on $K$ based at $x$ and the assignment $x\mapsto
\varphi _{x}$ is continuous in the compact--open topology.  It follows
rather easily that, for each $x\in K$ and each loop $\sigma _{x}$ on
$K$ based at $x$, the composed loop $\sigma _{x}\varphi _{x}\sigma
_{x}^{-1}$ is base point homotopic to $\varphi _{x}$.  This can only
be true if $\varphi _{x}$ is homotopically trivial, $\forall\,x\in K$.
Using this fact for each $x\in\bd K$, one constructs the desired
homotopy of $\varphi $ in $\homeo_{\,0} (K). $ We remark that the
group $\diff_{0} (K)$ is also known to be homotopically
trivial~\cite{E-S} as is the group of piecewise linear
homeomorphisms~\cite{scott_homeo}.
\end{rem}

Because of this lemma, we may write $\cf$ for $\cz$ and $\df$ for
$\dz$, where $f$ is the first return homeomorphism induced by $\LL$ on
a depth one leaf. We can also use $\zf$ to denote the isotopy class of
the core laminations corresponding to $f$. Here, the isotopies should
preserve each leaf of $\FF$.

\begin{lemma}\label{conj}
Let $f\co L_{0}\ra L_{0}$ be the first return homeomorphism induced on a
depth one leaf $L_{0}$ of $\FF$ by a transverse, $1$--dimensional
foliation $\LL$.  If $g\co L_{0}\ra L_{0}$ is a homeomorphism isotopic to
the identity, then $\cf=\CC_{gfg^{-1}}$.
\end{lemma}

Indeed, in standard fashion, the isotopy $g_{t}$, $g_{0}=g$ and
$g_{1}=\id$, induces an isotopy of $\LL$ to $\LL'$, leaving $\FF$
leafwise invariant, such that $\LL'$ induces first return map
$gfg^{-1}$ on $L_{0}$.

The set $Z=\ZZ\cap L_{0}$ is exactly the set of points which never
cluster at ends of $L_{0}$ under forward or backward iteration of the
monodromy $f$.  Assume that the dynamical system $(Z ,f)$ admits a
Markov partition $\{R_{1},\dots,R_{n}\}$ (in particular, these are
imbedded rectangles in $L_{0}$ that cover $Z $ and have disjoint
interiors) and let $(\Sigma _{A},\sigma _{A})$ be the associated
symbolic dynamical system.  Here, an $n\times n$ incidence matrix
$A=[a_{ij}]$ of 0's and 1's determines a closed subset $$\Sigma
_{A}\sseq\{1,2,\dots,n\}^{\Z},$$ a sequence $\iota
=\{i_{k}\}_{k=-\infty}^{\infty}$ being an element of $\Sigma _{A}$ if
and only if $a_{i_{k}i_{k+1}}=1$, $\forall\,k$. This is a compact,
metrizable, totally disconnected space and the shift map $\sigma _{A}$
is a homeomorphism.  In the usual scheme, there is a semiconjugacy
$$\varphi \co (\Sigma _{A},\sigma _{A})\ra(Z ,f)$$ defined by
$$\varphi (\iota )=\varphi (\cdots,i_{-1},i_{0},i_{1},\cdots) =
\bigcap_{k=-\infty}^{\infty}f^{-k}(R_{i_{k}})=R_{\iota }.$$  This
assumes that the infinite intersection $R_{\iota }$ of rectangles
$f^{-k}(R_{i_{k}})$ degenerates to a singleton, but we are going to
allow this set to be either a singleton, a nondegenerate arc, or a
nondegenerate rectangle.  We will still require that $$Z
=\bigcup_{\iota \in \Sigma _{A} }R_{\iota }, $$ but each symbol
sequence $$\iota =(\cdots,i_{-1},i_{0},i_{1},\cdots) \in \Sigma _{A}$$
will represent all the points in $R_{\iota }$.  Remark that a boundary
point of $R_{\iota }$ might be represented by distinct sequences in
$\Sigma _{A}$.

The closed orbits $\Gamma $ of $\Varphi _{t}$ determine periodic
orbits of $f$ in $Z $, hence correspond to periodic orbits of
$\sigma _{A}$. A point $\iota \in \Sigma _{A}$ has periodic $\sigma
_{A}$--orbit if and only if $\iota $ itself breaks down into a
bi-infinite sequence of a repeated finite string
$i_{0},\dots,i_{q-1}$, called a \textit{period} of $\sigma _{A}$. In
this case, $\iota $ is called a \textit{periodic point}. Given a periodic
point $\iota $, the Brouwer fixed point theorem implies that there is
at least one corresponding periodic $f$--orbit
$\{x,f(x),\dots,f^{q}(x)=x\}$, $x\in R_{\iota }$, and a corresponding
closed leaf $\Gamma _{\iota }=\{\Varphi _{t}(x)\}_{t\in\R}$ of $\ZZ$.

Let $\iota =\{i_{k}\}_{k=-\infty}^{\infty}\in \Sigma _{A}$ and suppose
that $i_{q}=i_{0}$ for some $q>0$. Let $x\in R_{\iota }$.  Then there
is a corresponding singular cycle $\Gamma _{q}$ formed from the
orbit segment $\gamma _{q}=\{\Varphi _{t}(x)\}_{0\le t\le q}$ and an
arc $\tau \ss R_{i_{0}}$ from $\Varphi _{q}(x)=f^{q}(x) \text{ to }
x$.  Also, since $i_{q}=i_{0}$, there is a periodic element $\iota
'\in
\Sigma _{A} $ with period $i_{0},\dots,i_{q-1}$ and a corresponding
closed leaf $\Gamma _{\iota '}=\Gamma '$ of $\ZZ$.

\begin{lemma}\label{periodic} 
The singular cycle $\Gamma _{q}$ and closed leaf $\,\Gamma '$, obtained as
above, are homologous. In particular, the homology class of $\Gamma _{q}$
depends only on the periodic element $\iota '$.
\end{lemma}

\proof
The loop $\Gamma '$ is the orbit segment $\{\Varphi _{t}(x')\}_{0\le
t\le q}$, for a periodic point $$x'\in R_{i_{0}}\cap
f^{-1}(R_{i_{1}})\cap\dots\cap f^{-q}(R_{i_{q}})=R'.$$ Remark that
$x\in R'$ also.  Let $\tau '$ be an arc in the subrectangle $R'\ss
R_{i_{0}}$ from $x$ to $x'$ and set $\tau ''=f^{q}(\tau ')$, an arc in
$f^{q}(R')$ from $f^{q}(x)$ to $x'$.  Since $i_{q}=i_{0}$,
$f^{q}(R')\ss R_{i_{0}}$ and the cycle $\tau +\tau '-\tau ''$ in the
rectangle $R_{i_{0}}$ is homologous to 0.  That is, we can replace the
cycle $\Gamma _{q}=
\gamma _{q}+ \tau $ by the homologous cycle $\gamma _{q}-\tau '+\tau
''$.  Finally, a homology between this cycle and $\Gamma '$ is
given by the map $$H\co [0,1]\times [0,q]\ra M,$$  defined by
parametrizing $\tau '$ on $[0,1]$ and setting
$$H(s,t)=\Varphi _{t}(\tau '(s)).\eqno{\qed}$$

If no proper, cyclicly consecutive substring of a $\sigma
_{A}$--period $i_{0},\dots,i_{q-1}$ also occurs as a period, we say
that the period is \textit{minimal}. It is elementary that there are
only finitely many minimal periods.  Those closed leaves $\Gamma $ of
$\ZZ$ that correspond to minimal periods in the symbolic system will
be called minimal loops in $\ZZ$.  The following is an easy
consequence of Lemma~\ref{periodic}.

\begin{cor}\label{minloops}
Every closed leaf $\Gamma $ of $\ZZ$ is homologous in $M$ to a linear
combination of the minimal loops in $\ZZ$ with non-negative integer
coefficients. Furthermore, every homology direction can be arbitrarily well
approximated by positive multiples of closed leaves of $\ZZ$.
\end{cor}

This corollary and Lemmas~\ref{span} and~\ref{periodic}
give the following important result.

\begin{theorem}\label{vertices} Suppose that the dynamical system
$(Z ,f)$ admits a Markov partition.  Then the cone $\cf\ss
H_{1}(M)$ is the convex hull of finitely many rays through classes
$[\Gamma_{i} ]$, $1\le i\le r$, where the structure cycles $\Gamma
_{i}$ are minimal loops in $\ZZ$.  Consequently, the dual cone $\df$
is polyhedral and both $\cf\text{ and }\df$ depend only on the
symbolic dynamics.
\end{theorem}


\section{Pseudo-Anosov endperiodic maps}\label{psA}

We continue with the hypotheses and notation of the preceding section.
Fix a noncompact leaf $L$ of $\FF$ and let $f\co L\ra L$ be the first
return map defined by a transverse 1--dimensional foliation $\LL_{f}$.
It is standard that $f$ is an \textit{endperiodic}
homeomorphism~\cite{fe:endp}.  Here, we use the well understood
structure theory of depth one leaves, writing $$ L= K\cup U_{+}\cup
U_{-}, $$ where $K$ is a compact, connected subsurface, called the
\textit{core} of $L$, $U_{\pm}$ falls into a disjoint union of finitely
many closed neighborhoods of isolated ends of $L$ and $K$ meets
$U_{\pm}$ only along common boundary components.  The set $U_{+}$ is
called the neighborhood of \textit{attracting ends} and has the
property that $f (U_{+})\ss U_{+}$.  The neighborhood $U_{-}$ of {\it
repelling ends} has the property that $U_{-}\ss f (U_{-})$. The core
is not unique since it can always be made larger by adjoining a
suitable piece of $U_{\pm}$. Set $K\cap U_{\pm}=\bd_{\pm}K$. While
$f^{n} (\bd_{+}K)$ suffers only bounded distortion as $n\ra\infty$, it
generally becomes unboundedly distorted as $n\ra-\infty$, the
situation being reversed for $f^{n} (\bd_{-}K)$.

The attracting ends $\{e^{i}_{+}\}_{i=1}^{r}$ of $L$ are permuted by
$f$, as are the repelling ends $\{e^{j}_{-}\}_{j=1}^{k}$.  The set of
cycles of these permutations corresponds one-to-one to the set of
components of $\tb M$.  Let
$\{e^{i_{1}}_{+},e^{i_{2}}_{+},\dots,e^{i_{q}}_{+}\}$ be a cycle
corresponding to the tangential boundary component $F$ and let
$U^{i_{1}}_{+}, U^{i_{2}}_{+},\dots,U^{i_{q}}_{+}$ be the
corresponding components of $U_{+}$ which are neighborhoods of these
ends. One can choose the above data so that \begin{align*} f
(U_{+}^{i_{j}}) &= U_{+}^{i_{j+1}},\quad 1\le j<q,\\ f (U_{+}^{i_{q}})
&\ss U_{+}^{i_{1}}.
\end{align*} There is a fundamental domain $F' \ss U_{+}^{i_{1}}$ for
the action of the semigroup $\{f^{n}\}_{n=0}^{\infty}$ on the union of
these neighborhoods. This domain $F' $ is homeomorphic to a manifold
obtained by cutting $F$ along the juncture
(cf~\cite[pages~3--4]{cc:smth2}) and it meets $K$ in a union
of common boundary components. A similar assertion holds for the
repelling ends of $L$, the semigroup being
$\{f^{-n}\}_{n=0}^{\infty}$.  The way in which the manifolds $f^{n}
(F' )$ link together along parts of their boundaries can be
surprisingly complicated.

As is well known, the depth one foliated manifold $(M,\FF)$ can be
recovered, up to foliated homeomorphism, from the endperiodic map
$f$. Indeed, the open, fibered manifold $(M_{0},\FF_{0}=\FF|M_{0})$ is
obtained, up to homeomorphism, by suspension of the homeomorphism $f$,
while the completion $$M=M_{0}\cup \tb M$$ is determined by the
endperiodic structure.  For more details see, for example,
\cite[Lemma~2.3]{cc:surg}. We further remark that the depth one
foliated manifold is homeomorphic (indeed, isotopic) to one in which
$\FF$ is smooth, even at the boundary, so we assume smoothness.

  The isotopy class $m (f)$ of $f$ (also called the \textit{mapping
class} of $f$) is completely determined by the depth one foliation
$\FF$ and, in turn, $m (f)$ determines the fibered manifold
$(M_{0},\FF_{0})$. The transverse foliation $\LL_{f}$ is not well
defined by $f$, although the isotopy class of its core lamination
$\ZZ_{f}$ is well defined (Lemma~\ref{1st:ret}). This isotopy class
varies, however, as $f$ is varied through endperiodic elements of $m
(f)$, so the cones $\cf$ generally change as $f$ is so varied.  We
want to choose $f$ so that these cones are as ``small'' as possible.
That is, we want the dual cone, $\df$ to be as large as possible.  The
tool for this is some unpublished work of Handel and Miller
(see~\cite{fe:endp}) which generalizes the Nielsen--Thurston
classification of homeomorphisms of compact surfaces~\cite{bca}. In
order to state this, some terminology is in order.

Denote by $L^{c}$ the compactification of $L$ obtained by adjoining
its ends.  By a \textit{properly imbedded line} in $L$, we mean a
topological imbedding $$\sigma\co [-\infty,\infty]\ra L^{c} ,$$ where
$\{\sigma (\pm\infty)\}$ is a pair of ends of $L$ and $\sigma
(-\infty,+\infty)\ss L$. This will be distinguished from a {\it
properly imbedded arc} which is an imbedding $$\sigma \co [-1,1]\ra
L,$$ where $\{\sigma (\pm1)\}=\sigma [-1,1]\cap\bd L$.  A {\it
peripheral curve} in $L$ is either a closed curve isotopic to a
component of $\bd L$ or a properly imbedded line, isotopic (with
endpoints fixed) to the endpoint compactification of a noncompact
component of $\bd L$.  A \textit{proper homotopy} between properly
imbedded lines or arcs is a homotopy that fixes the endpoints. If an
end $e$ of $L$ has a neighborhood that is homeomorphic either to
$S^{1}\times [0,\infty)$ or $[0,1]\times [0,\infty)$, then $e$ will be
called a \textit{trivial end.}

\begin{defn}
A closed, essential, nonperipheral curve $\gamma \ss L$ is a closed
reducing curve if, for a sufficiently large integer $n>0$, $f^{n}
(\gamma )\ss U_{+}$ and $f^{-n} (\gamma )\ss U_{-}$. A properly
imbedded, nonperipheral line $\sigma $ is a reducing line if one
endpoint is an attracting end, the other a repelling end, and $\sigma
$ is periodic under $f$ up to a proper homotopy. A periodic curve is a
closed, nonperipheral curve which is periodic under $f$ up to
homotopy.
\end{defn}

\begin{defn}\label{irred} The endperiodic map $f\co L\ra L$ is periodic
(or trivial) if every orbit $\{f^{n} (x)\}_{n=-\infty}^{\infty}$ has
points in $U_{+}$ and points in $U_{-}$.  The endperiodic map is
irreducible if  no end of $L$ is trivial and there
are no reducing lines, reducing curves, nor periodic
curves. Otherwise, $f$ is reducible.
\end{defn}

\begin{lemma}\label{compl_red}
If no component of $\tb M$ is an annulus or a torus, if $M$ is
completely reduced {\em (Definition~\ref{compl:red})}, and if $L$ is a
noncompact leaf of a taut, depth one foliation of $M$, then every
endperiodic homeomorphism $f\co L\ra L$ that occurs as the first return
map for a transverse foliation $\LL_{f}$ is irreducible. The
endperiodic homeomorphism is periodic if and only if $M\cong S\times
I$ and $\LL_{f}$ is a product $I$--bundle over $S$.
\end{lemma}

The elementary proof is left to the reader. The point of this lemma is
that our ongoing hypotheses (Theorem~\ref{ongoing}) imply that the
endperiodic monodromy is irreducible and nonperiodic.

 We can assume that the monodromy $f$ and associated transverse
foliation $\LL_{f}$ are smooth. Since no component of $\tb M$ is an
annulus or a torus, we will put a smooth Riemannian metric on $M$ such
that all the leaves are hyperbolic, chosen so that, in $U_{\pm}$, it
is the lift of the metric on $\tb M$ via the projection $U_{\pm}\ra\tb
M$ along the leaves of $\LL_{f}$~\cite{cc:prop2}.  We can demand that
$\bd L$ be geodesic and choose the juncture in each boundary leaf to
be geodesic. In particular, $\bd K$ will consist of geodesic arcs
and/or loops and $f$ will be an isometry in the ends.

\begin{theorem}[Handel--Miller]\label{HM}
Let $f\co L\ra L$ be an irreducible, endperiodic, nonperiodic monodromy
diffeomorphism of a hyperbolic leaf as above.  Then there is a pair of
mutually transverse geodesic laminations $\Lambda _{\pm}$ of $L$ and
an endperiodic homeomorphism $h\in m (f)$ such that
\begin{enumerate}
\item The components of $h^{n} (\bd_{\pm}K)$ are geodesics,
$\forall\,n\in\Z$, and
\begin{enumerate}
\item[{\em(a)}]  $\Lambda _{+}\ss K\cup U_{+}$ is the limit of the geodesic
laminations $h^{n} (\bd_{-}K)$ as $n\ra\infty$; \item[{\em(b)}] $\Lambda
_{-}\ss K\cup U_{-}$ is the limit of $h^{n} (\bd_{+}K)$ as
$n\ra-\infty$;
\end{enumerate}
 \item $K\cap(\Lambda _{+}\cup \Lambda _{-})$ weakly binds $K$ in the
sense that complementary regions not meeting $\bd K$ are simply
connected;
\item $h|\bd L=f|\bd L$; \item If $Z\ss K$ is the invariant set
of $h$, the dynamical system $(Z,h|Z)$ admits a Markov partition.
\end{enumerate}  Finally, $m (f)$ together with the hyperbolic
structure on $L$ uniquely determines $\Lambda _{\pm}$ and
$h|\Lambda _{+}\cap \Lambda _{-}$.  
\end{theorem}

\begin{rem}
In (a), the assertion that
$$
\lim_{n\ra\infty}h^{n} (\bd_{-}K) = \Lambda _{+}
$$ means that the leaves $\ell$ of $\Lambda _{+}$ are exactly the
curves in $L$ with the following property.  For each compact subarc
$J\ss\ell$, there is a sequence of subarcs $J_{k}\ss h^k (\bd_{-}K) $
converging uniformly to $J$ as $k\uparrow\infty$. Equivalently,
$\Lambda _{+}$ is the frontier of the open set
$$\UU_{-}=\bigcup_{n\in\Z}h^{n} (U_{-}).$$ This is the weak
interpretation of~(a).

For the strong interpretation of (a), let $\wt{L}$ be the the
universal cover of $L$ and let $\ol{L}$ be the closure of $\wt{L}$ in
the closed Poincar\'e disk $D\cup S_{\infty}$ (where $S_{\infty}$ is
the circle at infinity and $D$ is the open unit disk).  Either $\bd
L=\0$ and $\ol{L}=D\cup S_{\infty}$, or $\wt{L}$ has geodesic boundary
in $ D$ and $\ol{L}$ is the union of $\wt{L}$ and a Cantor set $C\ss
S_{\infty}$.  Since everything in sight is a geodesic, the weak
interpretation of (a) implies that the lifts $\wt{\ell}\ss\wt{L}$ of
leaves of $\Lambda _{+}$ are exactly the uniform limits in the {\it
Euclidean} metric of sequences $\wt{\sigma }_{k}$ of suitable lifts
of components of $h^{k} (\bd_{-}K)$. When $\wt{L}=D$, this means that
the endpoints in $S_{\infty}$ of (the completion of) $\wt{\sigma
}_{k}$ converge to the endpoints of $\wt{\ell}$.  In general, the
endpoints of $\wt{\sigma }_{k}$ may lie on geodesic boundary
components of $\wt{L}$, but they still converge in the Euclidean
metric to the endpoints of $\wt{\ell}$ in $C\ss S_{\infty}$.  Similar
remarks apply to~(b).
\end{rem}

We will say that $h\in m (f)$ is a pseudo-Anosov, endperiodic
automorphism of $L$.  Remark that $h$ preserves the laminations,
expanding the leaves of the unstable lamination $\Lambda _{+}$ and
contracting those of the stable lamination $\Lambda _{-}$. Unlike the
compact case, there may not be projectively invariant measures of full
support.

Since $\bigcup_{n=1}^{\infty}h^{n} (\bd_{-}K)$ consists of disjoint
geodesics and $K$ is a compact surface, the intersection of this set
with $K$ consists of compact, properly imbedded geodesic arcs which
fall into a finite number of isotopy classes rel~$\bd_{+}K$.  The
components of $\Lambda _{+}\cap K$ all belong to these isotopy
classes and every leaf of $\Lambda _{+}$ meets
$\bd_{+}K$, hence every leaf meets $h^{n} (\bd_{+}K)$, $n\ge1$.
Evidently, none of these leaves can meet $\bd_{-}K$.  Similar remarks
hold for the stable lamination $\Lambda _{-}$.

\begin{cor}\label{deep}
Every leaf $\ell$
of $\Lambda _{+}$ contains points of $Z$ and points going to
infinity in both directions in $\ell$ which converge to an attracting
end of $L$, but no points arbitrarily near repelling ends. The
analogous assertions hold for $\Lambda _{-}$, with the roles of
attracting and repelling ends interchanged.
\end{cor}

The points of $\Lambda _{+}\cap \Lambda _{-}$ are contained in the
core $K$ and remain there under all forward and backward iterations of
$h$.  Generally, this intersection is not the entire invariant set $Z$,
which may have nonempty interior~\cite[Proposition~2.12]{fe:endp}.

The lamination $\Lambda _{+}$ can be augmented to an $h$--invariant
geodesic lamination $\Gamma_{+}$ by adding on all $h^{n} (\bd_{-}K)$,
$-\infty<n<\infty$, and a similar augmented lamination $\Gamma_{-}$ is
obtained by adding all $h^{n} (\bd_{+}K)$ to $\Lambda _{-}$. Once the
choice of $K$ has been fixed, these augmented laminations are uniquely
determined by $m (f)$ and the hyperbolic metric, they are mutually
transverse and they weakly bind $L$.  The automorphism $h$ is unique
on $\Gamma _{+}\cap
\Gamma _{-}$ and may be extended continuously in any convenient way on the 
complementary arcs of this set in $\Gamma_{+}\cup
\Gamma_{-}$ and on the components of the complement of $\Gamma _{+}\cup
\Gamma_{-}$ in $L$. It is not clear that this can be done so
that $h$ will be smooth, even though the original monodromy $f$ was
smooth.

The foliation $\LL_{h}$, chosen to produce the pseudo-Anosov first
return map $h$, is a bit problematic at $\tb M$.  It happens that
there will always be distinct $h$--orbits in $\Lambda _{+}\cup \Lambda
_{-}$ which cluster in $M$ at the same point of $\tb M$
(see~\cite{fe:endp}). Each leaf of $\LL_{h}$ which comes close enough
to $\tb M$ actually limits on a unique point of $\tb M$, but the fact
that distinct leaves can limit on the same point of $\tb M$ implies
that the continuous extension of $\LL_{h}$ to $M $ will not be a
foliation.  We could remedy this by modifying $\LL_{h}$ in arbitrarily
small neighborhoods of $\tb M$ without affecting the core lamination
$\ZZ_{h}$, thereby allowing the isotopy arguments of the previous
section to be carried out.  Since the conclusions of these arguments
concern only the core lamination $\ZZ_{h}$, they remain true without
making such modifications.

By this remark and in order to keep the full force of pseudo-Anosov
monodromy, we agree that $\LL_{h}$ be defined only in $M_{0}$,
extending to $M$ as a singular foliation.

\begin{lemma}\label{2lam}
There are mutually transverse, smooth leaved, $2$--dimensional
laminations $\Theta _{\pm }$ of $M_{0}$ which intersect the leaves of
$\FF_{0}$ transversely in the augmented geodesic laminations
$\Gamma_{\pm}$ for the pseudo-Anosov monodromy of those leaves.
\end{lemma}

\begin{proof} For simplicity of exposition, we consider the case that
$\bd L=\0$.  This makes the universal cover of each leaf the full open
unit disk $D$. Standard modifications of the following argument prove
the general case.  Let $\Varphi _{t}$ be the leaf preserving flow in
$M_{0}$ with flow lines the leaves of $\LL_{f}|M_{0}$. Since $\LL_{f}$
is smooth, so is this flow.  Fix a leaf $L_{0}$ of $\FF|M_{0}$, set
$L_{t}=\Varphi _{t} (L_{0})$ and consider the open, $\FF$--saturated
set $V_{\epsilon }=\bigcup _{-\epsilon <t<\epsilon }L_{t}$. Here
$\epsilon >0$ is chosen small enough that $L_{\epsilon }\ne
L_{-\epsilon }$. The universal cover of $V_{\epsilon }$ is of the form
$D\times (-\epsilon ,\epsilon) $, where the interval fibers are flow
lines of the local flow $\wt{\Varphi }_{t}$ obtained by lifting
$\Varphi _{t}$, $-\epsilon <t<\epsilon $.

Write $D_{t}$ for $D\times
\{t\}$ and remark that the lifted metric gives a hyperbolic metric
$\gamma_{t} $ on $D_{t}$, $-\epsilon <t<\epsilon $.  Projection along
the interval fibers is smooth, but is not an isometry of these
metrics. However, for $\epsilon >0$ sufficiently small, projection
distorts the metrics only in a uniformly small way in the $\CII$
topology (indeed, this is true for $\Varphi _{t}$ downstairs, which
only fails to be an isometry on a compact neighborhood of $K$).  In
particular, each geodesic $\sigma_{0} $ in $D_{0}$ is carried to a
curve $\sigma _{t}$ in $D_{t}$ with uniformly small geodesic
curvature, $-\epsilon <t<\epsilon $.  An important property of
hyperbolic geometry is that any curve with geodesic curvature $<1$
remains uniformly near a unique geodesic, hence has well defined ends
in the circle at infinity. In fact, any pencil of asymptotic geodesics
limiting on the same point at infinity projects to a pencil of
``almost geodesics'' limiting on a common point at infinity. It
follows that projections along the interval fibers extend to well
defined projections on the circles at infinity. Thus, the completion $\MM$
of $D\times (-\epsilon ,\epsilon )$ by the circles at infinity of each
$D_{t}$ is canonically identified with $\ol{D}\times (-\epsilon
,\epsilon )$, making $\MM$ a smooth manifold.

Suppose that $\sigma _{0}$ is a geodesic covering a leaf $\ell_{0}$ of
$\Gamma _{\pm}$ in $L_{0}$ and consider the smooth surface
$S=\bigcup_{-\epsilon <t<\epsilon }\wt{\Varphi }_{t} (\sigma
_{0})$. If $z_{\pm\infty}$ are the endpoints of $\sigma _{0}$ at
infinity, then $S$ is continuously extended to $\ol{D}\times
(-\epsilon ,\epsilon ) $ by adjoining $\{z_{\pm\infty}\}\times
(-\epsilon ,\epsilon )$. Fix some third arc $\{w\}\times (-\epsilon
,\epsilon )$ at infinity.  If one uses the product metric on $D\times
(-\epsilon ,\epsilon )$ instead of the lifted metric, each leaf
$D_{t}$ has metric $\gamma _{0}$ and $\sigma _{t}=S\cap D_{t}$ is a
geodesic, $-\epsilon <t<\epsilon $.  But there is a unique leafwise
uniformizing diffeomorphism $\psi \co D\times (-\epsilon ,\epsilon )\ra
D\times (-\epsilon ,\epsilon )$ which fixes pointwise the three
transverse arcs at infinity and carries $\gamma _{0}$ to $\gamma
_{t}$, carrying each $\sigma _{t}$ to a $\gamma _{t}$--geodesic.  Each
of these geodesics in $D_{t}$ covers the leaf $\ell_{t}$ of $\Gamma
_{\pm}$ in $L_{t}$ which is in the homotopy class of $\Varphi _{t}
(\ell)$.  Since $\psi (S)$ is smooth, these leaves $\ell_{t}$ fit
together to form a smooth piece of surface transverse to $\FF$.
Extending this local construction in the obvious manner, we sweep out
the general leaf of the desired lamination $\Theta _{\pm}$.
\end{proof}

There is a certain amount of freedom in the choice of $h$ and
$\LL_{h}$. In order to use the Schwartzmann--Sullivan theory as in the
previous section, we note that Lemma~\ref{2lam} implies the following.

\begin{cor}\label{CO}
The foliation $\LL_{h}$, transverse to $\FF_{0}=\FF|M_{0}$ and
inducing pseudo--Anosov monodromy $h$ on a leaf $L$, can be chosen
to be integral to a $\CO$ vector field in $M_{0}$. 
\end{cor}

Indeed, $\Theta _{+}\cap \Theta _{-}$ is a smooth leaved,
1--dimensional lamination $X$ and we choose a continuous unit tangent
field to $X$.  This then extends, first to a continuous unit field on
$\Theta _{+}\cup \Theta _{-}$, thence to a field on $M_{0}$, which is
$C^{1}$ on each component of the complement of $X$ and is transverse
to the leaves of $\FF_{0}$.  Here we use the fact that $\Gamma_{+}\cup
\Gamma _{-}$ weakly binds $L$. The foliation integral to this field is
the desired realization of $\LL_{h}$.

We come to the key result of this section.

\begin{theorem}\label{also}
Let $h$ be a pseudo-Anosov first return map for $\FF$ and let $\FF'$
be another depth one foliation transverse to $\LL_{h}$, $f'$ the first
return map induced by $\LL_{h}$ on the typical noncompact leaf $L'$ of
$\FF'$.  Then there is a homeomorphism $g\co L'\ra L'$ which is isotopic
to the identity such that $h'=g\o f'\o g^{-1}\in m (f')$ is
pseudo-Anosov.
\end{theorem}

\begin{cor}\label{same cones}
Under the hypotheses of {\em Theorem~\ref{also}}, $\CC_{h}=\CC_{h'}$
and $\DD_{h}=\DD_{h'}$ and these cones are independent of the choice
of $h$ off of the set $\Lambda _{+}\cap
\Lambda _{-}$.
\end{cor}

Indeed, the first assertion is immediate by Lemma~\ref{conj} and the
second by Theorem~\ref{vertices}.

The proof of Theorem~\ref{also} will be broken down into a series of
lemmas. Again, for simplicity of exposition, we attend mainly to the
case that $\bd L=\0=\bd L'$ (equivalently, $\trb M=\0$), frequently
leaving it to the reader to adapt arguments to the general case.

To begin with, note that $\FF'$ will be transverse to
$\LL_{f}$ outside of a compact subset of $M_{0}$, so we fix a leafwise
hyperbolic metric for $\FF'$ such that, in neighborhoods of the ends
of depth one leaves, the metric is lifted from $\tb M$ by projection
along $\LL_{f}$.  

  Fix the noncompact leaves $L$ of $\FF$ and $L'$ of
$\FF'$. Projection along the leaves of $\LL_{h}$ defines local
homeomorphisms between $L\text{ and }L'$, but one cannot expect to
piece these together to a well defined covering map. Nonetheless, if
$s\co J\ra L$ is a curve, $t_{0}\in J$, one can project $s$ along
$\LL_{h}$ to a curve $s'\co J\ra L'$, this projection being completely
determined by the choice of $s' (t_{0})$.  Similarly, curves $s'$ on
$L'$ project to curves $s$ on $L$. The analysis
in~\cite[Section~4]{fe:endp} implies the following. 

\begin{lemma}\label{deep_in_ends} Deep  in the ends of $L\text{ and
}L'$ (that is, arbitrarily near $\tb M$ in the topology of
$M$), the projections along $\LL_{h}$ can be arbitrarily well
approximated by projections along $\LL_{f}$.  
\end{lemma}

It is sometimes helpful to lift the picture of this projection
operation to the universal cover $\wt{M}_{0}$ of $M_{0}$.  The
foliations $\FF_{0}\text{ and }\FF'_{0}$ lift to foliations
$\wt{\FF}_{0}\text{ and }\wt{\FF}'_{0}$ and, since the foliations
downstairs are taut, the lifted foliations have simply connected
leaves. The lift $\wt{\LL}_{h}$ of $\LL_{h}$ is transverse to both
lifted foliations.  Thus, we may view
$$
\wt{M}_{0} \cong \wt{L}\times \R \cong \wt{L}'\times \R,
$$ where $\wt{L}$ is a leaf of $\wt{\FF}_{0}$ covering $L$, $\wt{L}'$
is a leaf of $\wt{\FF}'_{0}$ covering $L'$, and the $\R$ factors are
leaves of $\wt{\LL}_{h}$.  Projection along the leaves of this
1--dimensional foliation carries $\wt{L}$ homeomorphically onto
$\wt{L}'$ and \textit{vice versa.}  This projection, restricted to a lift
of a curve in $L$ is a lift of a projection along $\LL_{h}$ as
described above.  This lifted picture of the projection is often
useful.  For instance, the following should be obvious.

\begin{lemma}\label{nullhomotopic}
Let $s\text{ and }s'$ be mutual projections along $\LL_{h}$.  Then
$s'$ is a nullhomotopic loop in $L'$ if and only if $s$ is a
nullhomotopic loop in $L$.
\end{lemma}

Let $\Lambda _{\pm}$ be the stable and unstable laminations preserved
by $h$ in the leaf $L$ of $\FF$. By taking all projections along
$\LL_{h}$ of the leaves of these laminations, we produce a pair of
laminations $\Lambda '_{\pm}$ on $L'$. It is clear that these will be
closed subsets of $L'$, but we cannot hope that they will be geodesic
laminations. We will show, however, that they have all the qualitative
properties that the pseudo-Anosov laminations should have.  It will
then be possible to construct a homeomorphism $g\co L'\ra L'$ which is
isotopic to the identity and simultaneously conjugates $\Lambda
'_{\pm}$ to the geodesic laminations of the Handel--Miller theory.
This will prove Theorem~\ref{also}.  Here again it will often be
useful to lift the data to the universal cover.

We begin with the analogue of Corollary~\ref{deep}.

\begin{lemma}\label{deep'}
Let $Z'\ss L'$ be the invariant set for $f'$.  Then every leaf $\ell$ of
$\Lambda' _{+}$ contains points of $Z'$ and points going to infinity in both
directions in $\ell$ which converge to an attracting end of $L'$, but no
points in a suitable neighborhood of the repelling ends. The analogous
assertions hold for $\Lambda '_{-}$, with the roles of attracting and
repelling ends interchanged.
\end{lemma}

\begin{proof}
Indeed, by Corollary~\ref{deep}, every leaf $\ell$ of $\Lambda _{+}$ contains
points of $Z$, hence meets $\ZZ_{h}$. It follows immediately that every leaf
of $\Lambda '_{+}$ meets $\ZZ_{h}$, hence contains points of $Z'$.
Similarly, no leaf of $\Lambda _{+}$ meets leaves of $\LL_{h}$ which approach
inwardly oriented components of $\tb M$, so the same holds for leaves of
$\Lambda '_{+}$.  That is, these leaves do not enter a periodic neighborhood
$U'_{-}$ of repelling ends.  Analogous remarks hold for the leaves of
$\Lambda '_{-}$.

Suppose that $\ell$ is a leaf of (say) $\Lambda '_{+}$ such that an end
$\alpha $ of $\ell$ has a neighborhood which does not enter a neighborhood of
the attracting ends of $L'$. The asymptote of $\alpha $ in $L'$ is therefore
a nonempty, compact sublamination $\lambda _{\alpha }$ of $\Lambda '_{+}$.
Evidently, the closure $\lambda _{n}$ in $L'$ of $$
\bigcup_{k=n}^{\infty} (f')^{-k} (\lambda _{\alpha })
$$ defines a descending nest $\{\lambda _{n}\}_{n=1}^{\infty}$ of nonempty,
compact sublaminations of $\Lambda '_{+}$ with intersection a compact,
nonempty, $f'$--invariant sublamination $\lambda $. Thus, $\lambda \ss Z'$,
implying that corresponding leaves of $\Lambda _{+}$ are contained in $Z$.
This contradicts Corollary~\ref{deep}.
\end{proof}

Let $\wt{\Lambda }_{\pm}$ and $\wt{\Lambda }'_{\pm}$ be the respective
lifts of our laminations to $\wt{L}$ and $\wt{L}'$.  Projection of
$\wt{L}'$ onto $\wt{L}$ along the leaves of $\wt{\LL}_{h}$ is a
homeomorphism carrying $\wt{\Lambda }'_{\pm}$ onto $\wt{\Lambda
}_{\pm}$.  Since the latter  are mutually transverse geodesic
laminations, the following is immediate.

\begin{lemma}\label{intersect}
Every leaf of $\wt{\Lambda }'_{+}$ meets every leaf of $\wt{\Lambda
}'_{-}$ in a single point.
\end{lemma}

It is somewhat touchier to prove that the leaves of $\wt{\Lambda
}'_{\pm}$ all limit nicely to endpoints at infinity.  There is no
guarantee that $\FF'$ is tangentially close to $\FF$, so there is no
way to guarantee that the laminations in $L'$ are leafwise close to
geodesic ones.

Consider an end $e$ of $L'$ and a neighborhood $W$ of $e$ in $L'$
which spirals on a leaf $F\ss\tb M$ over the juncture $N'\ss F$.  If
$\bd F=\0$, this juncture can be taken to be a simple closed geodesic.
In general, $\bd F\ne\0$ and $N'$ will be a disjoint union of simple
closed geodesics and/or properly imbedded geodesic arcs.  Similarly,
there is a juncture $N\ss F$ for the spiraling ends of the leaf $L$ of
$\FF$ which are asymptotic to $F$.  We work entirely in a normal
neighborhood $V=F\times [0,\epsilon )$ of $F$ in $M$ with normal
fibers subarcs of leaves of $\LL_{f}$, assuming this neighborhood to
be small enough that $\FF'| V$ is transverse to these normal fibers.
Thus, curves in $L\cap V$ can be projected both by $\LL_{f}$ and
$\LL_{h}$ to curves in $L'\cap V$ and \textit{vice versa.}

Let $\sigma $ be a component of $N'$ and fix a lift $\ol{\sigma
}$ of $\sigma $ to $L\cap V$ via projection along $\LL_{f}$. We
consider three cases. 

{\bf Case 1}\qua If $\sigma $ is a properly imbedded geodesic arc in $F$, so is
$\ol{\sigma }$ in $L$. Suitable projections $\sigma _{k}$ of $\ol{\sigma }$
into $L'$ along $\LL_{h}$ will be boundary components of complete
submanifolds $W_{k}\ss L'$ which form a fundamental neighborhood system of
$e$.  The reader can check that, in the upcoming arguments, this case can be
handled in close analogy with the next case.  In this way, we continue to
focus on the case in which the leaves have empty boundary.

{\bf Case 2}\qua If $\sigma $ is a closed geodesic with homological
intersection number $\sigma *N=0$, then $\ol{\sigma} $ will be a closed
geodesic in $L$.  Again, suitable projections $\sigma _{k}$ of $\ol{\sigma }$
into $L'$ along $\LL_{h}$ will be boundary components of complete
submanifolds $W_{k}\ss L'$ which form a fundamental neighborhood system of
$e$.  Note that each $\sigma _{k}$ has an $\LL_{h}$ projection to the
geodesic $\ol{\sigma }$ in $L$.

{\bf Case 3}\qua If $\sigma $ is a closed geodesic and $\sigma *N\ne0$,
then $\ol{\sigma }$ will be a half infinite geodesic.  We parametrize
it as $\ol{\sigma } (t)$, $0\le t<\infty$, remarking that, as
$t\uparrow\infty$, $\ol{\sigma } (t)$ moves arbitrarily far into an
end of $L$. Let $\ol{\sigma }_{\tau }$ denote the restriction of this
geodesic ray to $[\tau ,\infty)$.  For each $\tau \ge0$, suitable
$\LL_{f}$ projections $\sigma _{k}$ of $\ol{\sigma }_{\tau }$ are
closed geodesic boundary components of complete submanifolds $W_{k}\ss
L'$ which form a fundamental neighborhood system of $e$.  For suitably
large values of $\tau $ and large enough values of $k$, a suitable
projection $\sigma ' _{k}$ of $\ol{\sigma }_{\tau }$ into $L'$ along
$\LL_{h}$ will wrap around $\sigma _{k}$, staying as close to $\sigma
_{k}$ as desired. This assertion follows easily from
Lemma~\ref{deep_in_ends}. We emphasize that $\sigma '_{k}$ has an
$\LL_{h}$ projection to the geodesic ray $\ol{\sigma }_{\tau }$ in
$L$. \label{case3}

In Cases~2 and~3, let $\wt{\sigma}_{k} $ be a lift of $\sigma _{k}$ to
the universal cover, remarking that this lift limits on two distinct
points $a_{k}^{\pm}\in S_{\infty}$.  Suppose the lifts have been
chosen for all $k\ge K$ ($K>0$ fixed) so that there is a nested
sequence $$J_{K}\supseteq J_{K+1}\supseteq\cdots\supseteq
J_{k}\supseteq\cdots$$ of arcs in $S_{\infty}$ with $\bd
J_{k}=\{a_{k}^{+},a_{k}^{-}\}$. Thus, the simple closed curves
$\wt{\sigma }_{k}\cup J_{k}$ cut off a nested family $\{D_{k}\}_{k\ge
K}$ of closed disks in the closed Poincar\'e disk.

\begin{lemma}\label{contracts}
The nested sequence $\{D_{k}\}_{k\ge K}$ has intersection a singleton
$\{x\}$ in $S_{\infty}$.
\end{lemma}

\begin{proof} In both cases, the sequence
$\{\sigma _{k}\}_{k\ge K}$ converges to $e$. Thus, for every compact
subset $X\ss L'$, only finitely many terms of the sequence meet
$X$. The sequence of lifts $\{\wt{\sigma }_{k}\}_{k\ge K}$ must also
have the property that only finitely many terms meet any given compact
set $Y\ss\wt{L}'$.  In Case~3, $\sigma _{k}$ is a closed geodesic and,
in Case~2, the closed geodesic $\sigma '_{k}$ freely homotopic to
$\sigma _{k}$ stays uniformly close to $\sigma _{k}$ (again one uses
Lemma~\ref{deep_in_ends}).  In either case, the lifts of these
geodesics are complete geodesics in the hyperbolic plane $\wt{L}'$
having endpoints $a_{k}^{\pm}$ in the circle at infinity. This
sequence must also have the property that only finitely many terms
meet any given compact subset of $\wt{L}'$ and the assertion follows.
\end{proof}

In Case~1, The properly imbedded arcs $\sigma _{k}$ are homotopic
(with endpoints fixed) to properly imbedded geodesic arcs and the
reader can formulate and prove the appropriate analogue of
Lemma~\ref{contracts}. 

\begin{lemma}\label{not2}
If $\epsilon >0$, then there is an integer $K$ so large that, if $k\ge
K$ and $N_{k,\epsilon }$ is the normal neighborhood of hyperbolic 
radius $\epsilon $ of a lift $\wt{\sigma}_{k} $ of $\sigma _{k}$ to
$\wt{L}'$, then no leaf of $\wt{\Lambda }'_{\pm}$ can properly cross
$N_{k,\epsilon }$ twice.
\end{lemma}

\begin{proof}
Case~2 is the easier one.  In this case we can take $K=1$.  If
$\wt{\ell}$ is a leaf of $\wt{\Lambda }'_{\pm}$ which properly crosses
$\wt{\sigma} _{k}$ twice, there is a closed loop in $\wt{L}'$ of the
form $\wt{\tau} +\wt{\sigma} $, where $\wt{\tau }$ is an arc in
$\wt{\ell}$ and $\wt{\sigma }$ an arc in $\wt{\sigma }_{k}$.  The
covering map carries this to a nullhomotopic loop $\tau +\sigma $ in
$L'$, where $\tau $ is a subarc of a leaf $\ell$ of $\Lambda '_{\pm}$
and $\sigma $ is a path running around $\sigma _{k}$.  A suitable
projection along $\LL_{h}$ carries this to a nullhomotopic loop in $L$
(Lemma~\ref{nullhomotopic}) running around the closed geodesic
$\ol{\sigma }$ and along a subarc of a geodesic leaf of $\Lambda
_{\pm}$.  Evidently, these two geodesic segments do not coincide,
contradicting a well known property of hyperbolic surfaces.

In Case~3, one chooses $\tau >0$ and $K >0$ so large that, if $k\ge K$, a
suitable $\LL_{h}$ projection $\sigma '_{k}$ of the geodesic ray
$\ol{\sigma}_{\tau }$ stays inside the $\epsilon $--neighborhood of the
closed geodesic $\sigma _{k}$.  In the universal cover $\wt{L}'$, suitable
lifts $\wt{\sigma }'_{k}$ lie in $ N_{\epsilon ,k}$.  If $g$ is the
hyperbolic transformation with axis the complete geodesic $\wt{\sigma }_{k}$,
then $g^{n} (\wt{\sigma }'_{k})\ss N_{\epsilon ,k}$, $\forall\,n\in\Z$.  We
suppose that $\wt{\ell}$ is a leaf of $\wt{\Lambda }'_{\pm}$ properly
crossing $N_{\epsilon ,k}$ twice and note that, for suitable integers $n$,
$\wt{\ell}$ properly crosses $g^{n}(\wt{\sigma }'_{k})$ in two points.  Again
we obtain a nullhomotopic loop $\tau +\sigma $ in $L'$, where $\tau $ is a
subarc of a leaf $\ell$ of $\Lambda '_{\pm}$ and $\sigma $ is a path running
along $\sigma '_{k}$. A suitable $\LL_{h}$ projection carries this to a
nullhomotopic loop in $L$ made up of two distinct geodesic segments, one a
subarc of a leaf of $\Lambda _{\pm}$ and one a subarc of the geodesic ray
$\ol{\sigma }_{\tau }$.  This is the desired contradiction.
\end{proof}

\begin{cor}\label{infty}
Each leaf of $\wt{\Lambda }'_{\pm}$ limits on two distinct points in
the circle at infinity.
\end{cor}

\begin{figure}[ht!]
\vglue -0.7in
\begin{center}
\begin{picture}(300,150)(,-10)

\epsfxsize=300pt
\epsffile{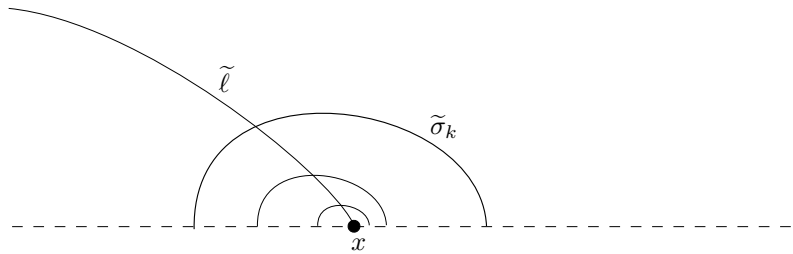}

\put (-170,-5){\small ${x}$}
\put (-220,55){\small ${\wt{\ell}}$}
\put (-140,39){\small ${\wt{\sigma }_{k}}$}

\end{picture}
\caption{Each end of $\wt{\ell}$ limits to a point at infinity in
$\wt{L}'$}\label{trapped} 
\end{center}
\end{figure}

\begin{proof}
Let $\wt{\ell} (t)$ cover the parametrized leaf $\ell (t)$ of $\Lambda
'_{\pm}$, $-\infty<t<\infty$ in $\wt{L}'$.  By Lem\-ma~\ref{deep'},
there is an end $e$ of $L'$ and a sequence $t_{k}\uparrow\infty$ such
that $\{\ell (t_{k}\}_{k=1}^{\infty}$ clusters at $e$ in $L'$. We can
assume that $\ell (t_{k})\in \bd W_{k}$ as above. Indeed, passing to a
subsequence, if necessary, and choosing the component $\sigma $ of
$N'$ appropriately, we can assume that $\ell (t_{k})\in\sigma _{k}$ as
in the above discussion.  Furthermore, as a little thought shows, we
can assume \textit{wlog} that, for infinitely many values of $k$, the
segment $\{\ell (t)\mid t_{k-1}<t<t_{k+1}\}$ does not meet $\sigma
_{k-1}$.  Thus, for all $k\ge1$, let $\wt{\sigma }_{k}$ be the lift of
$\sigma _{k}$ passing through the point $\wt{\ell} (t_{k})$. In the
case on which we are focusing, this lift will be an imbedded arc in
$D$ limiting on two distinct points $a_{k}^{\pm}\in S_{\infty}$.  In
Lemma~\ref{not2}, choose $\epsilon >0$ so small that $N_{k,\epsilon }$
cannot meet the lifts $\wt{\sigma }_{k\pm1}$.  Thus, for suitable
large values of $k$, there is no parameter value $t>t_{k-1}$ at which
$\wt{\ell}$ again meets $\wt{\sigma }_{k-1}$. This traps $\wt{\ell}
(t)$, $t>t_{k-1}$, in the closed disk $D_{k-1}$ of
Lemma~\ref{contracts}, so $$ \lim_{t\ra\infty}\wt{\ell} (t)=x
$$
 (see Figure~\ref{trapped}, where we use the upper half plane
model). Similarly the other end of $\wt{\ell}$ limits on a well
defined point $y$ in the circle at infinity.  If $x=y$, $\wt{\ell}$
would properly cross some $N_{k,\epsilon }$ twice, again contradicting
Lemma~\ref{not2}. \end{proof}

Let $K'\ss L'$ be a choice of core and $U'_{\pm}$ the corresponding
neighborhoods of attracting and repelling ends.  Let $\UU_{\pm }\ss L$
be the sets defined in the remark after Theorem~\ref{HM} and set $$
\UU'_{\pm} = \bigcup_{n\in\Z} (f')^{n} (U'_{\pm}).
$$ Then $\UU'_{-}$ is the set of points of intersection of $L'$ with
the leaves of $\LL_{h}$ which limit on inwardly oriented components of
$\tb M$ and $\UU_{-}$ is the corresponding set in $L$.  The sets
$\UU'_{+}$ and $\UU_{+}$ have an analogous description.  Since
$\Lambda _{\pm}$ is the set theoretic boundary in $L$ of $\UU_{\pm}$,
the following is immediate.

\begin{lemma}\label{set_th_bd}
The laminations $\Lambda '_{\pm}$ are the respective set theoretic
boundaries in $L'$ of $\UU'_{\pm}$.
\end{lemma}

Thus, \begin{align*}
\Lambda '_{+}&=\lim_{n\ra\infty} (f')^{n} (\bd_{-}K')\\ \Lambda
'_{-}&= \lim_{n\ra\infty} (f')^{-n} (\bd_{+}K')
\end{align*} in the weak sense.  It will be necessary to prove
convergence in the strong sense described in the remark following
Theorem~\ref{HM}.

Denote by $C$ the set of points at infinity of $\wt{L}'$.  We are
focusing on the case that $C=S_{\infty}$, but alternatively, it is a
Cantor set.  The following fact is well known for the universal
covering of a compact, hyperbolic surface, but is not generally true
for the noncompact case.  The proof for endperiodic surfaces was
communicated to us by S Fenley.

\begin{lemma}\label{min_at_infty} Let $G$ be the group of covering
transformations on $\wt{L}'$.  For each point $x\in\wt{L}'$, the orbit
$G (x)$ accumulates in $\wt{L}'\cup C$ exactly on the set $C$.  In
particular, $C$ is minimal under the induced action of $G$.
\end{lemma}

\begin{proof}
Let $x\in\wt{L}'$ and $z\in C$. It is clear that $G (x)$ cannot
accumulate in $\wt{L}'$ and we will show that $G (x)$ contains a
sequence converging to $z$ in the Euclidean metric.  Let $\gamma $
denote the unique unit speed geodesic ray in $\wt{L}'$ from $x$ to
$z$.  This projects to a geodesic ray or loop $\rho $ in $L'$,
parametrized on $[0,\infty)$, and we consider two cases.

{\bf Case I}\qua For a suitable sequence $t_{k}\uparrow\infty$, the hyperbolic
distance between $\rho (0)$ and $\rho (t_{k})$ is bounded by a finite
positive constant $B$.  Let $\rho _{k}=\rho |[0,t_{k}]$ and choose a
path $\tau _{k}$ in $L'$ from $\rho (t_{k})$ to $\rho (0)$ and having
length at most $B$.  Let $g_{k}\in G$ correspond to the lift of the
loop $\rho _{k}+\tau _{k}$ to a path starting at $x$.  Thus, this path
ends at the point $g_{k} (x)$ and is written as $\gamma_{k}+\wt{\tau
}_{k}$, where $\gamma_{k}=\gamma |[0,t_{k}]$ and $\wt{\tau }_{k}$ is
the lift of $\tau _{k}$ starting at $\gamma (t_{k})$.  The points
$\gamma (t_{k})$ converge to $z$ in the Euclidean metric, while the
Euclidean length of $\wt{\tau }_{k}$ converges to $0$. It follows that
$$
\lim_{k\ra\infty}g_{k} (x)=z.
$$

{\bf Case II}\qua The hyperbolic distance between $\rho (0)$ and $\rho
(t)$ goes to infinity with $t$. We can choose the sequence
$t_{k}\uparrow \infty$ so that $\rho (t_{k})$ converges to an end
$e$. We set $\rho_{k}=\rho |[0,t_{k}]$ and $\gamma _{k}=\gamma
|[0,t_{k}]$.  We can choose a fundamental neighborhood system
$\{W_{k}\}_{k=1}^{\infty}$ of $e$ with $\bd W_{k}$ geodesic and assume
that $\rho (t_{k})\in \sigma _{k}$, a component of $\bd W_{k}$,
$k\ge1$.  The lifts $\wt{\sigma }_{k}$ passing through $\gamma
(t_{k})$ cut off a fundamental neighborhood system
$\{D_{k}\}_{k=1}^{\infty}$ of $z$ in $\wt{L}'\cup C$.  Since the end
$e$ is nontrivial, there is a properly imbedded arc $\tau _{k}$ in
$W_{k}$ with endpoints in $\sigma _{k}$ which cannot be deformed into
$\sigma _{k}$ while keeping the endpoints in $\sigma _{k}$.  We can
deform $\tau _{k}$, keeping the endpoints in $\sigma_{k} $, to a loop
(again denoted by $\tau _{k}$) based at $\rho (t_{k})$ and let
$g_{k}\in G$ correspond to the lift of $\rho _{k}+\tau _{k}-\rho _{k}$
starting at $x$. This lift has the form $\gamma _{k}+\wt{\tau
}_{k}-\wt{\rho }_{k}$, where $\wt{\tau }_{k}$ is the lift of $\tau
_{k}$ starting at $\gamma (t_{k})\in\wt{\sigma }_{k}$ and $-\wt{\rho
}_{k}$ is the lift of $-\rho _{k}$ starting at the terminal point of
$\wt{\tau }_{k}$. By the assumption on $\tau _{k}$ this terminal point
cannot lie on $\wt{\sigma }_{k}$ and it follows that $g_{k} (x)\in
D_{k}$.
\end{proof}

Consider a leaf $\ell$ of $\Lambda '_{\pm}$ and a lift $\wt{\ell}$ of
this leaf to $\wt{L}'$.  For definiteness, assume $\ell$ to be a leaf
of $\Lambda '_{+}$.  Since $(f')^{k} (\bd_{-}K')$ converges to
$\Lambda '_{+}$ in the weak sense, we can find a component $\sigma $
of a juncture arbitrarily deep in a repelling end and a sequence of
positive integers $n_{k}\uparrow\infty$ such that arbitrarily long
subarcs of $\sigma _{k}=(f')^{n_{k}} (\sigma )$ uniformly well
approximate arbitrarily long subarcs of $\ell$ as $k\uparrow\infty$.
This likewise holds for the lift $\wt{\ell}$ and for suitable choices
of the lifts $\wt{\sigma }_{k}$. Note that $\wt{\sigma }_{k}$ has well
defined endpoints $a_{k}^{\pm}$, either in $\bd \wt{L}'$ or in $C$ and
that Corollary~\ref{infty} guarantees that $\wt{\ell}$ has well
defined endpoints in $C$.  We can assume that there are well defined
limits $$ a^{\pm}=\lim_{k\ra\infty}a_{k}^{\pm}.
$$ In order to prove the strong convergence, we need only show that
$a^{\pm}\in C$ are the endpoints of $\wt{\ell}$.  We suppose not and
reach a contradiction.

\begin{figure}[ht!] 
\vglue -0.2in
\begin{center}
\begin{picture}(300,100)(,) 
\epsfxsize=300pt 
\epsffile{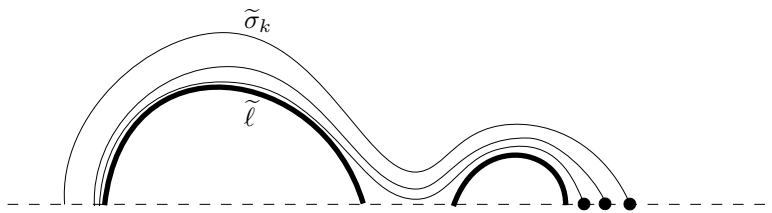}
\put (-210,33){\small$\wt{\ell}$}
\put (-210,70){\small$\wt{\sigma }_{k}$}
\end{picture} 
\caption{Nonhausdorff accumulation of $\{\wt{\sigma }_{k}\}_{k=1}^{\infty}$
on leaves of  
$\wt{\Lambda }_{+}'$}\label{nonhaus}
\end{center} 
\end{figure}

No generality is lost in assuming that the component $\sigma $ of
juncture is chosen as in one of the three cases discussed in the
remarks preceding Lemma~\ref{contracts}. There, $\sigma $ denoted a
component of juncture in a boundary leaf $F$, but here we fix a
homeomorphic lift deep in a repelling end of $L'$ and denote that lift
by $\sigma $.  By Lemma~\ref{min_at_infty} and the invariance of
$\wt{\Lambda }'_{+}$ under the covering group $G$, we conclude that
the endpoints of the leaves of this lamination are dense in $C$.
Since $\wt{\sigma }_{k}$ cannot meet these leaves, we easily conclude
that the sequence $\{\wt{\sigma }_{k}\}_{k=1}^{\infty}$ has
nonhausdorff accumulation on more than one leaf of $\wt{\Lambda
}'_{+}$ (actually, on infinitely many) as indicated in
Figure~\ref{nonhaus}. The following, therefore, completes the proof of
convergence in the strong sense. The idea for this proof was suggested
to us by S Fenley.

\begin{lemma}\label{haus}
The sequence $\wt{\sigma }_{k}$ converges uniformly to $\wt{\ell}$ in
the Euclidean metric on $D$.
\end{lemma}

\begin{proof}
We suppose nonhausdorff accumulation as in Figure~\ref{nonhaus} and deduce a
contradiction.  Let $\pi \co \wt{L}'\ra\wt{L}$ denote the projection along the
leaves of $\wt{\LL}_{h}$.  This is a homeomorphism.  It will preserve the
nonhausdorff accumulation, although we are no longer assured that $\pi
(\wt{\sigma }_{k})$ has well defined endpoints at infinity.  Since $\pi
(\wt{\ell})$ is a geodesic leaf of $\wt{\Lambda }_{+}$, it does have well
defined endpoints. Let $z$ be one of these endpoints, remarking that $\pi
(\wt{\sigma }_{k})$ passes arbitrarily near $z$ (in the Euclidean metric) as
$k\ra\infty$.

Let $\ol{\ell}$ denote the leaf of $\Lambda _{+}$ covered by $\pi
(\wt{\ell})$.  By Corollary~\ref{deep}, there is a sequence
$\{W_{k}\}_{k=1}^{\infty}$ of closed neighborhoods of an attracting end of
$L$ such that $\bd W_{k}$ is geodesic and such that $\ol{\ell} (t_{k})\in\bd
W_{k}$, where $\pi (\wt{\ell} (t_{k}))\ra z$ as $k\ra\infty$.  Thus, we can
fix a fundamental system $\{D_{k}\}_{k=1}^{\infty}$ of closed neighborhoods
of $z$ in $\wt{L}$, each bounded by a geodesic (or geodesic arc) $\wt{\tau}
_{k}$ which is a lift of the component $\tau _{k}$ of $\bd W_{k}$ which
passes through $\ol{\ell} (t_{k})$.

\begin{figure}[ht!]
\begin{center}
\vglue-0.3in
\begin{picture}(300,100)(,)

\epsfxsize=300pt
\epsffile{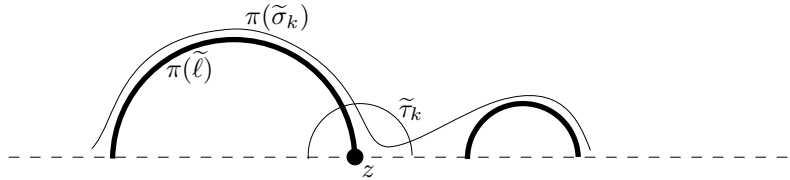}
\put (-210,55){\small$\pi (\wt{\sigma }_{k})$}
\put (-152,20){\small$\wt{\tau }_{k}$}
\put (-240,35){\small$\pi (\wt{\ell})$}
\put (-166,-3){\small$z$}

\end{picture}
\caption{$\pi
(\wt{\sigma }_{k})$ crosses the geodesic $\wt{\tau }_{k}$
twice}\label{nullhomotopy} 
\end{center}
\end{figure}

Since the curves $\pi (\wt{\sigma }_{k})$ also accumulate on other leaves of
$\wt{\Lambda }_{+}$, there will be large enough values of $k$ so that $\pi
(\wt{\sigma }_{k})$ crosses $\wt{\tau }_{k}$ twice
(Figure~\ref{nullhomotopy}).  Thus, we obtain a nullhomotopic loop in
$\wt{L}$ made up of a segment of $\pi (\wt{\sigma }_{k}) $ and a segment of
the geodesic $\wt{\tau }_{k}$.  If $p\co \wt{L}\ra L$ denotes the covering
projection, we obtain a nullhomotopic loop $\rho _{k}$ in $L$ made up of a
segment of the geodesic $\tau _{k}$ and of the curve $p (\pi(
\wt{\sigma }_{k}))$.

If $\sigma $ is as in Case~2, then there is an $\LL_{h}$ projection of
$\sigma $ to a geodesic loop $\ol{\sigma }$ deep in a repelling end of
$L$.  Since $f'$ is everywhere defined by projection along $\LL_{h}$,
every $\sigma _{k}$ has an $\LL_{h}$ projection to the closed geodesic
$\ol{\sigma }$.  Since $h$ is also defined by projection along
$\LL_{h}$, it follows that a suitable negative iterate $h^{-n_{k}}$
carries $p (\pi( \wt{\sigma }_{k}))$ onto the geodesic $\ol{\sigma }$.
Since $h$ preserves the extended geodesic lamination $\Gamma _{-}$,
$h^{-n_{k}} (\tau _{k})$ is also geodesic and $h^{-n_{k}} (\rho_{k} )$
is a nullhomotopic loop made up of two distinct geodesic segments, a
contradiction.  Case~1 is entirely similar.

If $\sigma $ is as in Case~3, one chooses a normal $\epsilon
$--neighborhood $N_{\epsilon } (\sigma )$ of the closed geodesic
$\sigma $ so that all iterates $N_{\epsilon,k }=(f')^{k} (N_{\epsilon
} (\sigma ))$ are pairwise disjoint.  The geodesic ray $\ol{\sigma
}_{\tau }\ss L$ in Case~3 has an $\LL_{h}$ projection contained
entirely in $N_{\epsilon } (\sigma )$. In the argument above, view the
lifts $\wt{N}_{\epsilon ,k}$ as slightly thickened versions of
$\wt{\sigma }_{k}$, also accumulating in nonhausdorff fashion on
multiple leaves of $\wt{\Lambda }^{'}$.  Modifications of the above
argument, analogous to the treatment of Case~3 in the proof of
Lemma~\ref{not2}, produce a nullhomotopic loop made up of a segment of
the geodesic $\ol{\sigma }_{\tau }$ and a segment of the geodesic
$h^{-n_{k}} (\tau _{k})$.
\end{proof}

Since the leaves of $\wt{\Lambda }'_{+}$ are exactly the uniform limits of
lifts of components of $(f')^{k} (\bd_{-}K')$ and all of these curves have
well defined (and distinct) endpoints in $\wt{L}'\cup C$, we can replace each
of these curves with the geodesic having the same endpoints.  The strong
convergence proven above implies corresponding convergence of endpoints,
hence strong convergence of the geodesics with these endpoints.  Thus, these
geodesics are exactly the lifts of the leaves of the extended geodesic
laminations $\Gamma '_{+}$ of the Handel--Miller theory.  Similar arguments
apply to the lamination $\Gamma '_{-}$.  In standard fashion, one constructs
a homeomorphism $\wt{g}\co \wt{L}'\ra\wt{L}'$, commuting with the covering
group, fixing $C$ pointwise, and carrying the laminations to geodesic ones.
This induces the homeomorphism $g\co L'\ra L'$ of Theorem~\ref{also}, completing
the proof of that theorem.

\begin{cor}\label{existence}
Let $(M,\FF)$ be a depth one, foliated, sutured manifold.  Then there is a
unique closed cone $\DD\ss H^{1} (M)$ which is convex, polyhedral,
contains the proper foliated ray $[\FF]$, has interior a union of foliated
rays, and is such that $\bd\DD$ contains no foliated rays.
\end{cor}

\begin{proof}
By Theorem~\ref{reduced}, Lemma~\ref{ann-tor} and Lemma~\ref{compl_red},
together with Theorem~\ref{HM}, we can suppose that the endperiodic monodromy
$h$ of a leaf $L$ of $\FF$ is pseudo-Anosov and consider the cone $\DD=\dh$.
By Theorem~\ref{vertices}, this cone is polyhedral. By Theorem~\ref{dual
cone}, every class in $\intr\dh$ is foliated and we only need to prove that
there is no foliated class in $\bd\dh$. The linear inequalities $\theta
_{i}\ge0$ defining $\dh$ are given by \textit{integral} homology classes
$\theta _{i}\in H_{1} (M;\Z)$, these being represented by closed loops, so
there is a foliated class $[\omega ]\in\bd\dh$ if and only if some such class
corresponds to a depth one foliation. Let $\FF'$ be a depth one foliation
represented by a class $[\omega' ]\in\bd\dh$ and let $h'\co L'\ra L'$ be a
pseudo-Anosov first return map for a noncompact leaf $L'$ of $\FF'$.  Choose
a class $[\omega '']\in\intr \dh$ on a ray through the integer lattice and as
close to $[\omega ']$ as desired.  Thus, we can assume that the closed,
nonsingular foliated form $\omega ''$ is sufficiently near $\omega '$ that
the depth one foliation $\FF''$ which it defines is also transverse to
$\LL_{h'}$. But Corollary~\ref{same cones} then implies that $\dh=\DD_{h'}$,
so $[\omega ']\in\intr \dh$, a contradiction.
\end{proof}

\section{Finiteness of the Foliation Cones.}

One can construct nonmaximal foliation cones in close analogy with
U.\ Oertel's construction of branched surfaces which carry norm
minimizing representatives of elements of $H^{1}
(M)$~\cite{Oer:hom}.  While it seems to be hard to use this
approach to prove convexity of the maximal cones, it does provide a
proof of finiteness.  Indeed, every maximal foliation cone is a union
of Oertel cones and there are only finitely many of these.  The proof
requires a few modifications of arguments in~\cite{Oer:hom} since it
will not be convenient to use norm minimizing surfaces.  The norm
minimizing hypothesis in Oertel's argument is used to prove a crucial
orientation property~\cite[pages~261--262]{Oer:hom} which, in our case,
will follow from the depth one hypothesis. We assume familiarity with
the construction of branched surfaces from surfaces in Haken normal
form as carried out, for example, in~\cite{F-O,Oer,Oer:hom}.

We fix a handlebody decomposition of the sutured manifold $M$,
requiring that $\bd M$ be covered by a union of $i$--handles, $0\le
i\le2$, $\tb M\cap\trb M$ being covered by a union of 0--handles and
1--handles. All constructions of surfaces in normal form and of
branched surfaces will be carried out relative to this handlebody
decomposition. 

We also fix normal neighborhoods $W=\tb M\times [0,2]$
and $V=\trb M\times [0,2]$, arranging that the handles meeting $\tb M$
lie in $W'=\tb M\times [0,1]$ and those meeting $\trb M$ lie in
$V'=\trb M\times [0,1]$.

Given a depth one foliation $\FF$, we will produce a branched surface
$B$ which, in a suitable sense, ``fully carries'' $\FF$. (For this, we
could invoke a theorem of Gabai~\cite[Theorem~4.11]{gab:kneser}, but
we will describe a more elementary proof.) 

The branched surface $B$ will be such that every choice $\mu =\{\mu
_{i}\}$ of strictly positive weights on the sectors $B_{i}$ of $B$,
satisfying the \textit{branch equations}~\cite[page 386]{Oer}, defines a
foliated class $[\mu] \in H^{1} (M)$. We call $\mu $ an
\textit{invariant measure} on $B$.

As in the cited references, one directly constructs the normal
neighborhood $N (B)\ss M$, foliated by compact, oriented intervals
$I$, $B$ being defined as the quotient of $N (B)$ obtained by
collapsing the interval fibers to points.  While $N (B)$ is not
exactly an interval bundle in the usual sense, the local models
pictured in the cited references make clear the sense in which it can
be viewed as an interval bundle over $B$.  Hereafter, we use the term
``interval bundle'' in this more general sense.

\begin{figure}[ht!]
\begin{center}
\begin{picture}(300,160)(-40,-10)

\epsfxsize=200pt
\epsffile{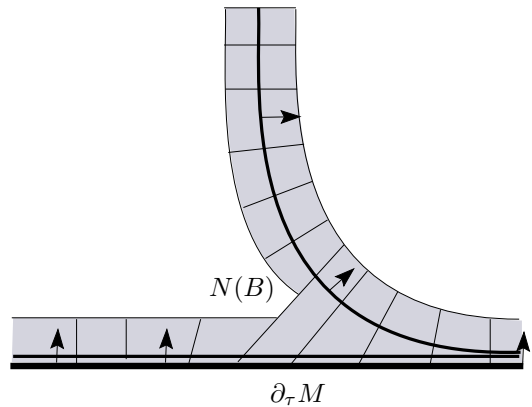}
\put (-100,-10){\small$\tb M$}
\put (-123,30){\small$N (B)$}

\end{picture}
\caption{The normal neighborhood $N (B)$ near $\tb M$\label{N(B)}}
\end{center}
\end{figure}

A major difference from the construction in~\cite{Oer:hom} is that
$\tb M$ will be part of $B$.  This implies that $B$ can fully carry
only noncompact surfaces $L$ which accumulate on $\tb M$ exactly as
does a depth one leaf (Figure~\ref{N(B)}).  This requires that the
corresponding measure $\mu $ take the value $\infty$ on the sectors of
$B$ contained in $\tb M$. Away from the tangential boundary, $L$ will
meet each interval fiber of $N (B)$ in only finitely many points.

The foliation $\FF$ will be transverse to the fibers of the normal
neighborhood $\tb M\times [0,\epsilon ]\ss W'$ for small enough
$\epsilon >0$, so an isotopy that flows into $M$ along the fibers of
$W$ allows us to assume that the foliation $\FF|W$ is transverse to
these fibers. Here, each noncompact leaf of $\FF|W$ falls into
finitely many disjoint, noncompact pieces, each spiraling in on one or
another component of $\tb M$.  It will also be convenient to choose
the fibers of $V=\trb M\times [0,2]$ to be tangent to the leaves.

Choose a depth one leaf $L$ of $\FF$. This submanifold, being a leaf
of a taut foliation, is incompressible in $M$. As
in~\cite[Section~4]{cc:isotopy}, fix a decomposition $L=C\cup E$,
where $C$ is a compact, connected submanifold of $L$ with boundary
and, possibly, corners, $E$ being the union of spiraling pieces.  We
can assume that $E$ is exactly the intersection of $L$ with the
handles that meet $\tb M$.  As in~\cite[pages~166--167]{cc:isotopy},
complete $C$ to a properly imbedded surface $F=C\cup A$ in $M$, where
$A$ is made up of annuli and/or rectangles which ``drop'' from
components of $\bd C$ to $\tb M$. Since $L$ is incompressible, so is
$F$.

We would like to perform an isotopy on $F$, putting it into Haken
normal form relative to the handlebody decomposition.  As the referee
pointed out, there is a problem applying standard theory here since we
cannot guarantee that $F$ is boundary incompressible.  However, the
very simple structure of $F$ in the neighborhood $U=W'\cup V'$ of $\bd
M$ makes it clear that there are no boundary compressing disks
contained entirely in $U$, so standard methods do permit an isotopy
that puts $A$ into Haken normal form relative to the handles meeting
$\bd M$. We view this as an ambient isotopy, simultaneously moving the
leaf $L$.  Now the incompressibility of $F$ allows us to extend this
to an ambient isotopy, supported in the complement of the interiors of
the handles meeting $\bd M$ and putting all of $F$, hence $A$, into
Haken normal form. The fact that $E=L\sm\intr C$ consists of spiraling
pieces near $\tb M$ now allows the entire leaf $L$ to be isotoped to
Haken normal form.

Let $N (L)\ss M$ be an imbedded normal neighborhood of $L$, this being
an oriented $J$--bundle over $L$, where the natural orientation of the
compact interval $J\ss\R$ agrees with the transverse orientation of
$L$. We can assume that $N (L)$ itself is in normal form.  That is, $N
(L)$ meets no 3--handles and, if $H$ is an $i$--handle, $1\le i\le 2$,
then each component of $H\cap N (L)$ is of the form $D\times J$, where
$D$ is one of finitely many \textit{disk types} in $H$. Let $N^{+}
(L)$ be the augmented $J$--bundle, incorporating a normal neighborhood
of $\tb M$ which lies in the union of $i$--handles meeting $\tb M$ and
is also in normal form. The orientation of the fibers respects the
transverse orientations of $L\text{ and }\tb M$ as indicated in
Figure~\ref{N(B)}.  Because $L$ is a depth one leaf, $M_{0}\sm \intr
N^{+} (L)$ will have a product structure $G\times I$, where $G$ is a
compact 2--manifold and $I$ a compact interval.

 In~\cite{Oer:hom}, the complexity of a compact, properly imbedded
surface in Haken normal form is defined to be the number of disks in
which the surface intersects 2--handles. Since $L$ meets the
2--handles along $\tb M$ in infinitely many disks, we modify this
definition, using only the 2--handles which do not meet the tangential
boundary. We can assume that, among all normal forms in its isotopy
class, $L$ has minimal complexity.  Following Oertel, we form
$\ol{N}^{+} (L)$ by adjoining to $N^{+} (L)$ all products $D\times
[-1,1]$, where $D\times
\{\pm1\}$ are adjacent, normally isotopic disks of $\bd N (L)\cap H$,
$H$ an $i$--handle with $0\le i\le2$. We can write $$
\ol{N}^{+} (L)= N^{+} (L)\cup Q_{1}\cup\cdots\cup Q_{r},
$$ where the $Q_{i}$'s are compact, connected, pairwise disjoint products.
As is standard, we break $\bd\ol{N}^{+} (L)$ into the horizontal part
$\bd_{h}\ol{N}^{+} (L)$ and the vertical part $\bd_{v}\ol{N}^{+} (L)$.
Again, $\ol{N}^{+} (L)$ is a compact $J$--bundle and the complement in
$M_{0}$ of its interior is a product (possibly not connected).

\begin{figure}[ht!]
\vglue -0.3in
\begin{center}
\begin{picture}(300,130)(,)

\epsfxsize=300pt
\epsffile{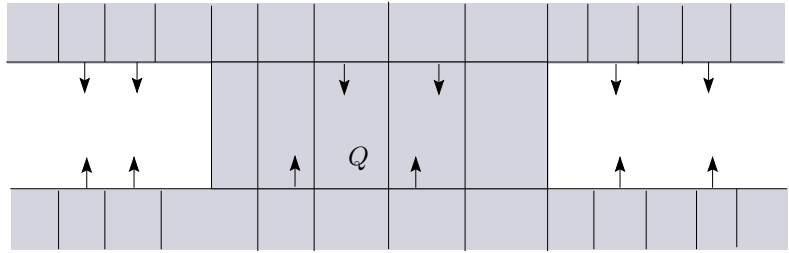}
\put (-170,35){\small$Q$}

\end{picture}
\caption{Opposed orientations on $\bd_{h} Q$}\label{opposed}
\end{center}
\end{figure}

Following Oertel~\cite{Oer,Oer:hom}, we cut $\ol{N}^{+} (L)$ along
$\bd N (L)\cap Q_{j}$, $1\le j\le r$, obtaining a compact interval
bundle $\wh{N}^{+} (L) $ which generally has many components. There is
a canonical immersion $\iota \co \wh{N}^{+} (L)\ra M$ which is an
imbedding on each component of $\wh{N}^{+} (L)$ and has image
$\ol{N}^{+} (L)$. We can identify the components of $\wh{N}^{+} (L)$
with their images under $\iota $, noting that one component of
$\wh{N}^{+} (L)$ is identified with $N^{+} (L)$.  One must be
concerned with the possibility that some component(s) $Q=Q_{j}$ may
have opposed orientations along the two components of $\bd_{h} Q$, as
indicated in Figure~\ref{opposed}.

\begin{lemma}\label{oriented}
There is a continuous choice of orientation on the interval fibers of
 $\ol{N}^{+} (L)$ which agrees with the orientation on the fibers of
 $N^{+} (L)$.
\end{lemma}

\begin{figure}[ht!]
\vglue -0.3in
\begin{center}
\begin{picture}(300,130)(,)

\epsfxsize=300pt
\epsffile{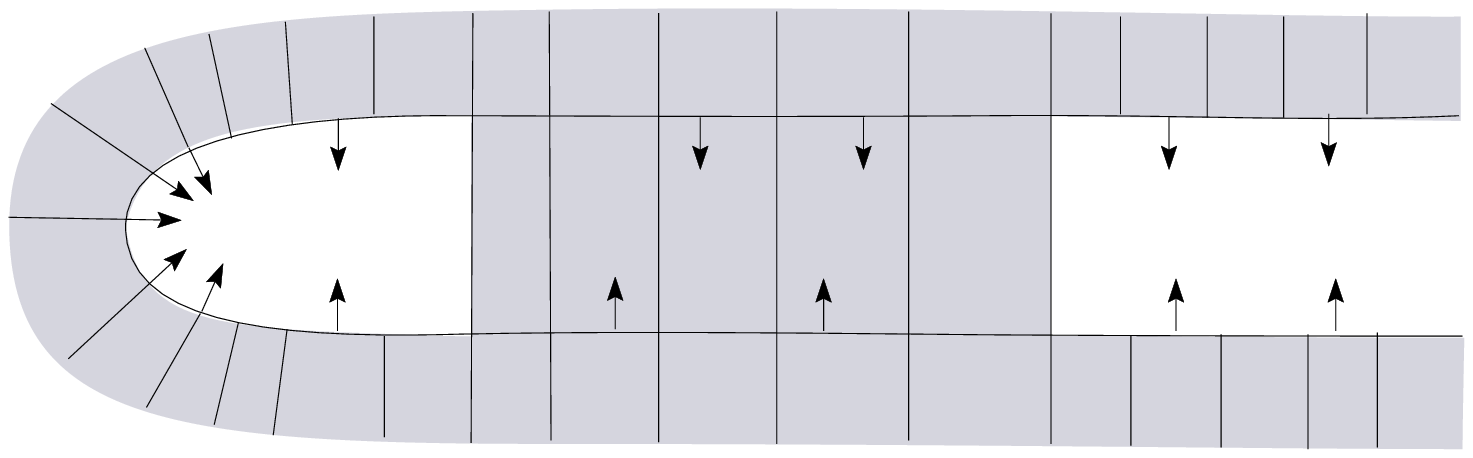}
\put (-160,45){\small$Q$}
\put (-240,45){\small$B$}

\end{picture}\nocolon
\caption{}\label{homotopy}
\end{center}
\end{figure}

\begin{figure}[ht!]
\vglue -0.3in
\begin{center}
\begin{picture}(300,100)(,)

\epsfxsize=300pt
\epsffile{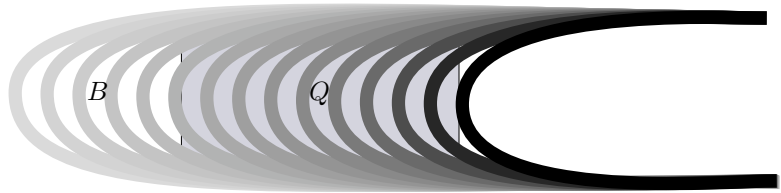}
\put (-184,38){\small$Q$}
\put (-268,38){\small$B$}

\end{picture}
\caption{A complexity reducing isotopy}\label{isotopy}
\end{center}
\end{figure}

\begin{proof}
We assume that some component $Q=Q_{j}$ of $\wh{N}^{+} (L)$ has
opposed orientations on $\bd_{h}Q$, as indicated in
Figure~\ref{opposed}, and deduce a contradiction.  Write $Q=R\times
[0,1]$ and suppose that $C\ss R$ is any closed loop (respectively,
properly imbedded arc).  Then $C\times [0,1]$ can be interpreted as a
proper homotopy of $C\times \{0\}$ to $C\times \{1\}$ in the
complement of the depth one leaf $L$ in $M_{0}$. Because of the
opposed orientations on $\bd_{h}Q$, the homotopy meets that leaf at
$C\times \{0,1\}$ \textit{on the same side of $L$}. Since $M_{0}\sm L$ is
a product, this homotopy can be compressed to a homotopy in $L$,
keeping $C\times \{0,1\}$ pointwise fixed.  Remark that $R\times
\{0\}$ and $R\times \{1\}$ are disjoint subsurfaces of $L$ and $C$ is
arbitrary, so the well understood structure of orientable surfaces
implies that these subsurfaces are annuli (respectively, disks) and
are separated by an annulus (respectively, a disk).  One concludes
that the picture is as in Figure~\ref{homotopy}, where $B\text{ and
}Q$ are both solid tori (respectively, 3--cells). An isotopy of $N
(L)$ across $B\cup Q$, as indicated in Figure~\ref{isotopy}, then
produces a normal form in the isotopy class of $L$ with strictly
smaller complexity.  This contradicts the fact that $L$ already has
minimal complexity.
\end{proof}

\begin{cor}\label{LL}
There is a smooth, oriented, one dimensional foliation $\LL$ of $M$ which is
transverse to $\tb M$, tangent to $\trb M$ and in $\ol{N}^{+} (L)$ coincides
with the interval fibration.
\end{cor}

Indeed, $P=M_{0}\sm\intr\ol{N}^{+} (L)$ is a product, so the corollary
follows directly from Lemma~\ref{oriented}. Note that $\LL|P$ is the
product fibration of $P$ by compact intervals.

The branched surface $B$ is the quotient of $\ol{N}^{+} (L)$ by the
interval fibers.  We write $\ol{N}^{+} (L)=N (B)$ and view $B$ as
imbedded in $N (B)$ transverse to the fibers in the standard way.  By
a small isotopy of $N (B)$ we can assume that the branch locus meets
itself transversely and without triple points.  

\begin{rem}
In~\cite{Oer:hom}, certain ``trivial'' components of $\wh{N}^{+} (L)$
are deleted, there\-by eliminating disks of contact and insuring that
$B$ is Reebless as well as incompressible.  For our purposes, there
seems to be no compelling reason to incorporate this step.
\end{rem}

If an imbedded surface $F\ss N (B)$ is transverse to the fibers, meets
every fiber and meets the fibers away from $\tb M$ only finitely
often, $F$ is said to be \textit{fully carried} by $B$.  Because $\tb
M$ can be viewed as part of $B$ with branch locus the junctures, any
surface fully carried by $B$ will be noncompact with ends spiraling in
on $\tb M$ over these junctures.  Such a surface defines an invariant
measure $\mu $ on $B$ which has value $\infty$ on sectors of $B$ in
$\tb M$ and positive integer values on other sectors. Conversely, an
invariant measure with these properties corresponds to such an
imbedded surface.  More generally, we allow the finite values of $\mu
$ to be positive real numbers, always requiring that the branch
equations be satisfied.  

An invariant measure $\mu $ well defines a homomorphism $$
[\mu ]\co \pi _{1} (M)\ra\R.
$$ Indeed, given a loop $\sigma $ in $M$, there is a free homotopy of
$\sigma $ to a loop $\sigma '$ which is transverse to $B$ and meets
$B$ only in the ineriors of sectors $B_{i}$.  If $x\in \sigma' \cap
B_{i}$, we assign the value $\pm\mu _{i}$ to $x$ according as the
orientation of $\sigma' $ at $x$ does or does not agree with the
transverse orientation of $B$ there.  The sum of these signed weights
is an invariant of the free homotopy class of $\sigma $ because of the
branch equations.  Thus, $[\mu ]\in H^{1} (M)$.

\begin{lemma}\label{foliated_class}
Each cohomology class $[\mu ]$, determined as above by a strictly
positive invariant measure $\mu $ on $B$, is a foliated class. As $\mu
$ varies over all such measures, these classes form a convex cone
$\DD_{B}$ in $H^{1} (M)$.
\end{lemma}

\begin{proof} 
Let $\ZZ$ be the core lamination of the foliation $\LL$ of
Corollary~\ref{LL} and let $\Varphi _{t}$ parametrize $\LL$ as a flow
(stationary exactly along $\tb M$). One can use this parametrization
to define the ``length'' of any compact subarc of a leaf of $\LL$,
remarking that there is a positive upper bound $m$ to the lengths of
the fibers of $N (B)$ that do not meet $\tb M$ and to the lengths of
the fibers of $P$.  Let $\ol{\Gamma }= \lim_{k\ra\infty}\Gamma
_{k}^{*}/\tau _{k}$ be a homology direction of $\ZZ$, where $\Gamma
_{k}^{*}/\tau _{k}$ are \textit{long, almost closed orbits} as defined
in the proof of Lemma~\ref{long_aco}. If $a>0$ is the minimum value taken by
$\mu $, one checks that the values of $[\mu ]$ on the long, almost
closed orbits cluster at a value no smaller than $a/2m$.  By
Lemma~\ref{span}, it follows that $[\mu ]$ is strictly positive on the
cone $\Cz$ of structure cycles. In particular, none of these structure
cycles bound, so $[\mu ]\in\intr\dz$ is a foliated class
(Theorem~\ref{dual cone}).  The fact these foliated classes form a
convex subcone $\DD_{B}\ss\dz$ is elementary.
\end{proof}

The branched surfaces $B$ obtained by our procedure are finite in
number, up to small isotopy.  Indeed, each $i$--handle determines
finitely many disk types, up to transverse isotopy, and contributes at
most one disk of each type to the construction of $B$.  Since
$\DD_{B}$ is a convex cone of foliated classes, it is contained in a
unique maximal foliation cone as in Theorem~\ref{existence}. This
guarantees that there are only finitely many of the latter.

\begin{theorem}
There are only finitely many maximal foliation cones.
\end{theorem}

The proof of Theorem~\ref{cones} is complete. The assertion that the
proper foliated ray $[\FF]$ determines $\FF$ up to isotopy was proven
in~\cite{cc:isotopy}.

\begin{rem}
The ``Oertel cone'' $\DD_{B}$ may not have full dimension, but at
least some of the branched surfaces produce cones of full
dimension. Indeed, the nonempty interior of each foliation cone is a
finite union of Oertel cones.  \end{rem}

\section{Computing examples}

If $M$ reduces to a sutured manifold that is completely disk
decomposable~\cite{ga0}, the foliation cones are easy to compute.
Furthermore, in many examples it is easy to compute the Thurston norm
for the 3-manifold $M$ using the foliation cones. In these examples,
each foliation cone is the union of cones over top dimensional faces
of the Thurston ball and the set $$\mathfrak{C} = \{ C= (-C_i)\cap C_j
\mid C_i,C_j\text{ are foliation cones}\}$$ is the set of cones over
the top dimensional faces of the Thurston ball. If $M^{*}$ denotes the
sutured manifold obtained by reversing the orientation of $\tb M$ (but
not of $M$) the cones $$\{-C_i \mid C_i\text{ is a foliation cone for
}M\}$$ become the foliation cones for $M^{*}$. The lattice points in
the cones $(-C_{i})\cap C_{j}$ correspond to foliations of both
$M\text{ and }M^{*}$. Thus, in the examples mentioned above, the
lattice points in the interior of each cone $C\in \mathfrak{C}$
correspond to foliations that can be spun both ways at $\tb M$. In
general it is very hard to compute the Thurston norm and these results
do not hold.

\begin{rems}
We routinely make the identification $$H_{2} (M,\bd M)=H^{1} (M).$$
When representing foliated classes by disks of a disk decomposition,
it is more natural to view the foliation cones in $H_{2} (M,\bd M)$.
When defining these cones by inequalities $\Gamma _{\iota }\ge0$,
where $\Gamma _{\iota }$ is (the homology class of) a loop in $\ZZ$ of
minimal period $\iota $, it is more natural to view the cones in
$H^{1} (M)$.

We compute the Thurston norm in a sutured manifold by doubling along
the sutures, computing the norm of the doubled class in the doubled
manifold, and dividing by $2$. Up to a factor $2$, this is the same as
Scharlemann's definition of the Thurston norm for a sutured
manifold~\cite[Definition~7.4]{scharle} (take $\beta =\0$ in
Scharlemann's definition).

In the following examples if $\{R_1,\ldots,R_n\}$ is a Markov
partition for $\ZZ\cap L$, we will let $a_{ij}=1$ iff $h(R_i)\cap
R_j\ne \emptyset$ (zero otherwise), and $A =
\left[a_{ij} \right]$. Then $$\Sigma_A = \{ (\dots i_k\,i_{k+1}\dots)
\mid a_{i_k 
i_{k+1}} = 1,\forall\,k\in\Z\}$$ are the allowable sequences (see
Section~\ref{structurecycles}). If $\iota \in \Sigma _{A}$ is
periodic, the transverse loop $\Gamma _{\iota }$ is well defined as a
homology class.

 If the sutured manifold $M$ is completely disk
decomposable, using disks $$D_1,\ldots,D_n\ss M,$$ we will
let $\alpha_i$ be the curves in $M$ defined, up to homology, by the
condition that the intersection products with the disks are
$$\alpha_i\cdot D_j = \delta_{ij},\quad 1\le i,j \le n.$$ We will use
the notation $\mathbf{e}_{i} = [D_i]\in H_2(M,\bd M)$, $0\le i\le n$.
(In many of the examples there will be a disk $D_0$ so that ${\mathbf{e}_0} + {\mathbf{e}_1} + \ldots + {\mathbf{e}_n} = 0$). Then $H_2(M,\bd M)
=\R^n$ is generated by ${\mathbf{e}_1},\ldots,{\mathbf{e}_n}$. The minimal
period loops $\Gamma _{\iota }$, thought of as classes in $H_{1} (M)$,
will be expressed in terms of the basis $\{\alpha _{1},\dots,\alpha
_{n}\}$.

In our examples, the sutured manifolds will be of the form $M_S(\kappa
)$ as in the introduction, where $\kappa $ is a knot or link.  They
will generally be denoted simply by $M$ unless some confusion is
possible. They will be pictured ``from the inside''.  That is, $M_{S}
(\kappa )$ will be the complement of the interior of the handlebody
that is drawn.  Our convention will be that $S_{+}$ is the part of
$\tb M$ oriented outwards from $M$ (away from the viewer), $S_{-}$ the
part oriented inwards.  A decomposing disk $D_{i}$ with positive sign
(often suppressed) has boundary oriented counterclockwise from the
viewer's perspective, the boundary of $-D_{i}$ being oriented
clockwise.
\end{rems}

\begin{figure}[ht!]
\begin{center}
\begin{picture}(300,150)(0,70)

\epsfxsize=330pt
\epsffile{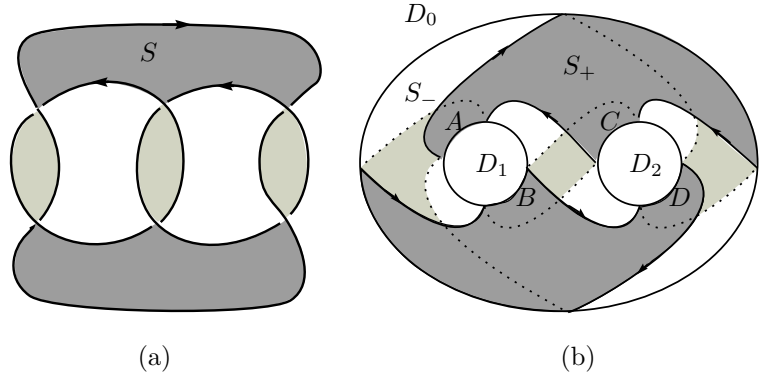}
\put (-265,70){\small (a)}
\put (-105,70){\small (b)}
\put (-265,185){\small$S$}
\put (-165,200){\small$D_{0}$}
\put (-165,170){\small$S_{-}$}
\put (-105,180){\small$S_{+}$}
\put (-150,159){\small$A$}
\put (-123,130){\small$B$}
\put (-138,143){\small$D_{1}$}
\put (-80,143){\small$D_{2}$}
\put (-91,159){\small$C$}
\put (-65,130){\small$D$}

\end{picture}
\caption{(a) A Seifert surface for $(2,2,2)$\quad (b) The sutured
manifold $M$ obtained from $(2,2,2)$}\label{(2,2,2)} 
\end{center}
\end{figure}

Finally, the following observation will often be used.

\begin{prop}\label{twice}
Suppose that $$ (M,\gamma ) \stackrel{D}{\leadsto} (M',\gamma ') $$ is
a disk decomposition by a nonseparating disk that meets the sutures
twice. Then  the regluing map $$ p\co M'\ra M $$
induces a surjection $$ p^{*}\co H^{1} (M)\ra H^{1} (M') $$ having
one-dimensional kernel generated by the Poincar\'e dual of
$D$. Furthermore, the foliation cones in $H^{1} (M)$ are exactly the
preimages of the foliation cones in $H^{1} (M')$.
\end{prop}

Indeed, an elementary Mayer--Vietoris argument, expressing $M$ as the
union of $M'$ and $D\times I$,  proves the assertion.  This is similar
to the procedure in Section~\ref{reducing}, but easier and more
natural for disk decomposable examples.

\begin{figure}[ht!]
\begin{center}
\begin{picture}(300,250)(10,-10)

\epsfxsize=320pt
\epsffile{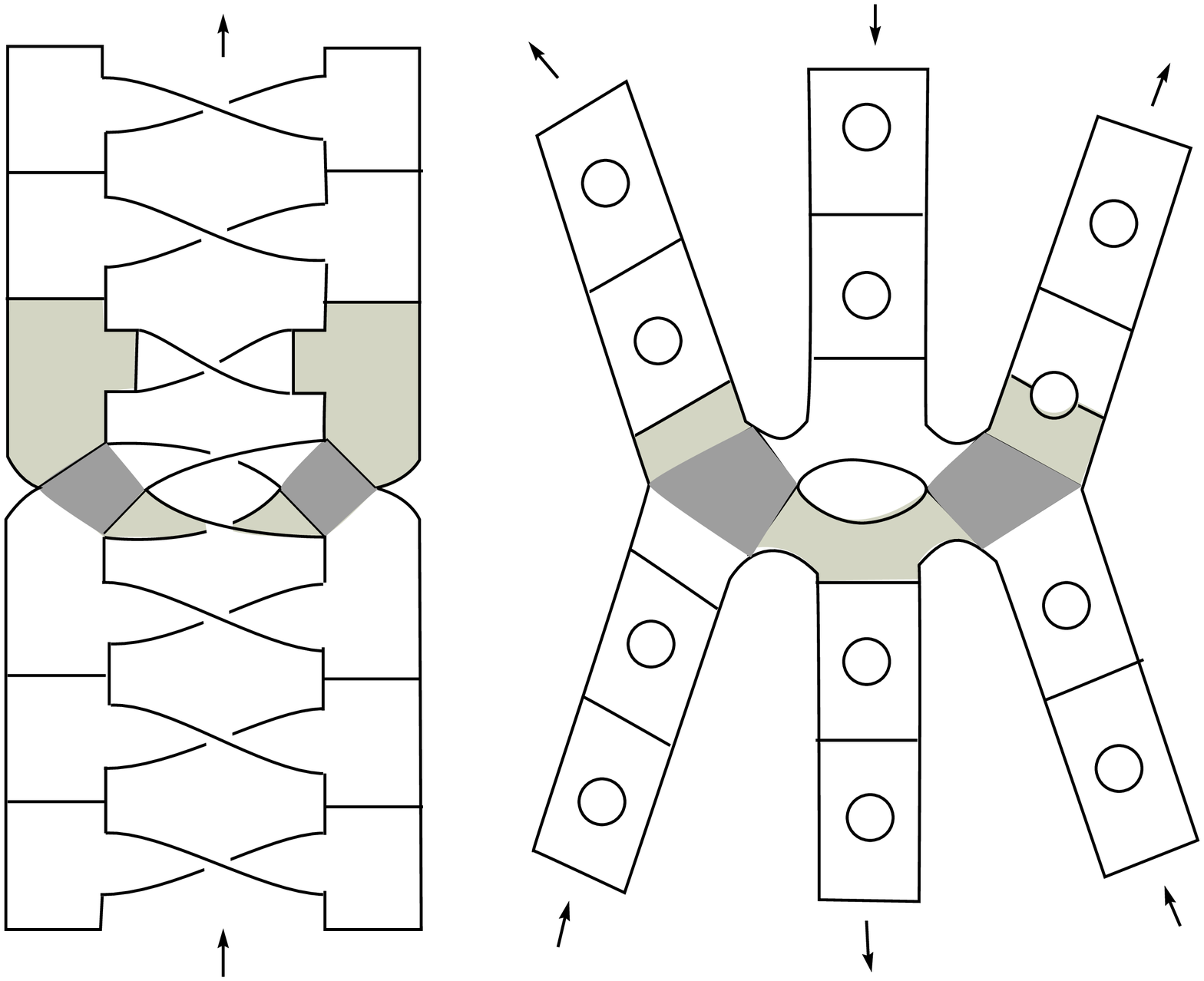}
\put (-106,0){\small (b)}
\put (-254.5,0){\small (a)}
\put (-220,130){\small$A$}
\put (-242,114){\small$B$}
\put (-272,114){\small$C$}
\put (-290,130){\small$D$}
\put (-226.5,171){\small$h (D)$}
\put (-298,171){\small$h (A)$}
\put (-233,149){\small$h (B)$}
\put (-293,149){\small$h (C)$}

\end{picture}
\caption{(a) The leaf L of Example~\ref{exone}\quad (b) The leaf L
of Example~\ref{extwo}}\label{leafL} 
\end{center}
\end{figure}

\begin{example}\label{exone}

Cutting the complement of the three component link $(2,2,2)$ apart
along the Seifert surface given in Figure~\ref{(2,2,2)}(a), gives the
sutured manifold $M$, everything outside the solid two-holed torus in
Figure~\ref{(2,2,2)}(b). We let $D_1\text{ and }D_2$ be the two disks
indicated in the figure and $D_0$ the outside disk.

\begin{figure}[ht!]
\vglue-0.3in
\begin{center}
\begin{picture}(300,130)(0,80)

\epsfxsize=320pt
\epsffile{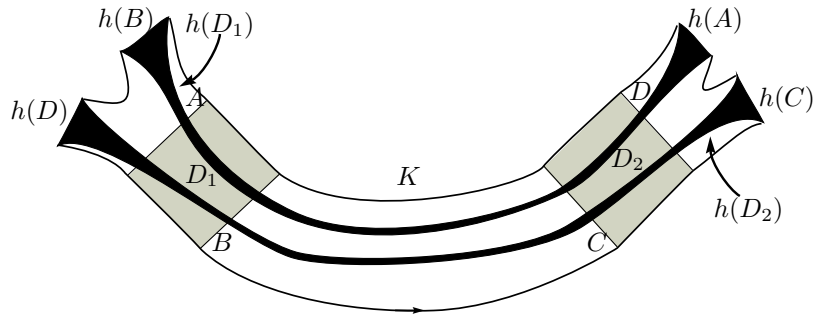}
\put (-170,125){\small$K$}
\put (-250,155){\small$A$}
\put (-240,100){\small$B$}
\put (-250,126){\small$D_{1}$}
\put (-250,183){\small$h (D_{1})$}
\put (-283,185){\small$h (B)$}
\put (-60,185){\small$h (A)$}
\put (-50,113){\small$h (D_{2})$}
\put (-33,155){\small$h (C)$}
\put (-316,150){\small$h (D)$}
\put (-82,157){\small$D$}
\put (-89,132){\small$D_{2}$}
\put (-98,100){\small$C$}

\end{picture}
\caption{The core $K$ of the leaf $L$ of
Example~\ref{exone}}\label{(2,2,2)core} 
\end{center}
\end{figure}

If one carries out the disk decomposition of $M$ using the disks
$+D_1,+D_2$ to obtain a foliation $\FF$ of $M$, one can piece together
a core $K$ of a leaf $L$ of $\FF$ as displayed in
Figure~\ref{(2,2,2)core}, the entire leaf $L$ being displayed in
Figure~\ref{leafL}(a) with the core shaded.  (It should be remarked
that this choice of $K$ is not entirely in accord with the conventions
of Sections~\ref{structurecycles} and~\ref{psA}, where the components
of $L\sm K$ are all unbounded.  The current choice, however, is
adequate to cover the invariant set $Z$ and is convenient for
illustrating the pseudo-Anosov dynamics.) In Figure~\ref{leafL}(a),
the endperiodic map $h$ moves points near the ends in the directions
indicated by the arrows, but also involves a lateral exchange as
indicated by the labels. With some thought, this can all be deduced
from Figure~\ref{(2,2,2)}. In Figure~\ref{(2,2,2)core}, the disks
$D_{1}\text{ and }D_{2}$ serve as the rectangles of a Markov
partition.  The intersections $h (D_{i})\cap D_{j}$ are also
displayed, showing that the incidence matrix is $$ A = \begin{bmatrix}
1&1\\ 1&1
\end{bmatrix}.$$
That is, in this case every sequence is an allowable sequence. The
loops corresponding to sequences of minimal period are \begin{align*}
\Gamma _{(\dots 111 \dots)}&= \alpha_1\\ \Gamma _{(\dots
222 \dots)} &= \alpha_2,
\end{align*}  so the foliation cone to which $\FF$ belongs  is
defined by the inequalities \begin{align*}
\{\alpha _{1}\ge 0\} &=\text{the right
half plane}\\ \{\alpha _{2}\ge 0\} &=\text{the upper half plane.}
\end{align*}
A base for this cone is the segment $[{\mathbf{e}_1},{\mathbf{e}_2}]$ (see
Figure~\ref{(2,2,2)cones}).

\end{example}

There is a simple symmetry of $M$ cyclically permuting
${\mathbf{e}_1},{\mathbf{e}_2},{\mathbf{e}_0}$. This symmetry takes
the foliation cone with base $[{\mathbf{e}_1},{\mathbf{e}_2}]$ to the
foliation cone with base $[{\mathbf{e}_2},{\mathbf{e}_0}]$ and then to
the foliation cone with base $[{\mathbf{e}_0},{\mathbf{e}_1}]$. This
gives Figure~\ref{(2,2,2)cones}, where the unit ball for the Thurston
norm on $M$ is also drawn with dashed line segments.

\begin{figure}[ht!]
\begin{center}
\begin{picture}(300,120)(10,240)

\epsfxsize=350pt
\epsffile{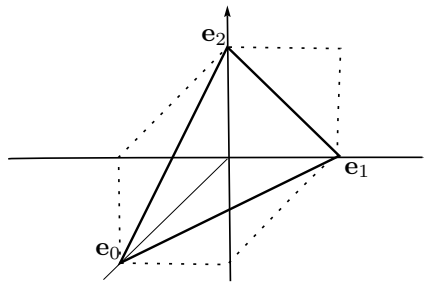}
\put (-224,259){\small$\mathbf{e}_{0}$}
\put (-130,290){\small$\mathbf{e}_{1}$}
\put (-184,340){\small$\mathbf{e}_{2}$}

\end{picture}
\caption{Foliation cones for $(2,2,2)$}\label{(2,2,2)cones}
\end{center}
\end{figure}

\begin{figure}[ht!]
\vglue-0.3in
\begin{center}
\begin{picture}(300,150)(0,70)

\epsfxsize=320pt
\epsffile{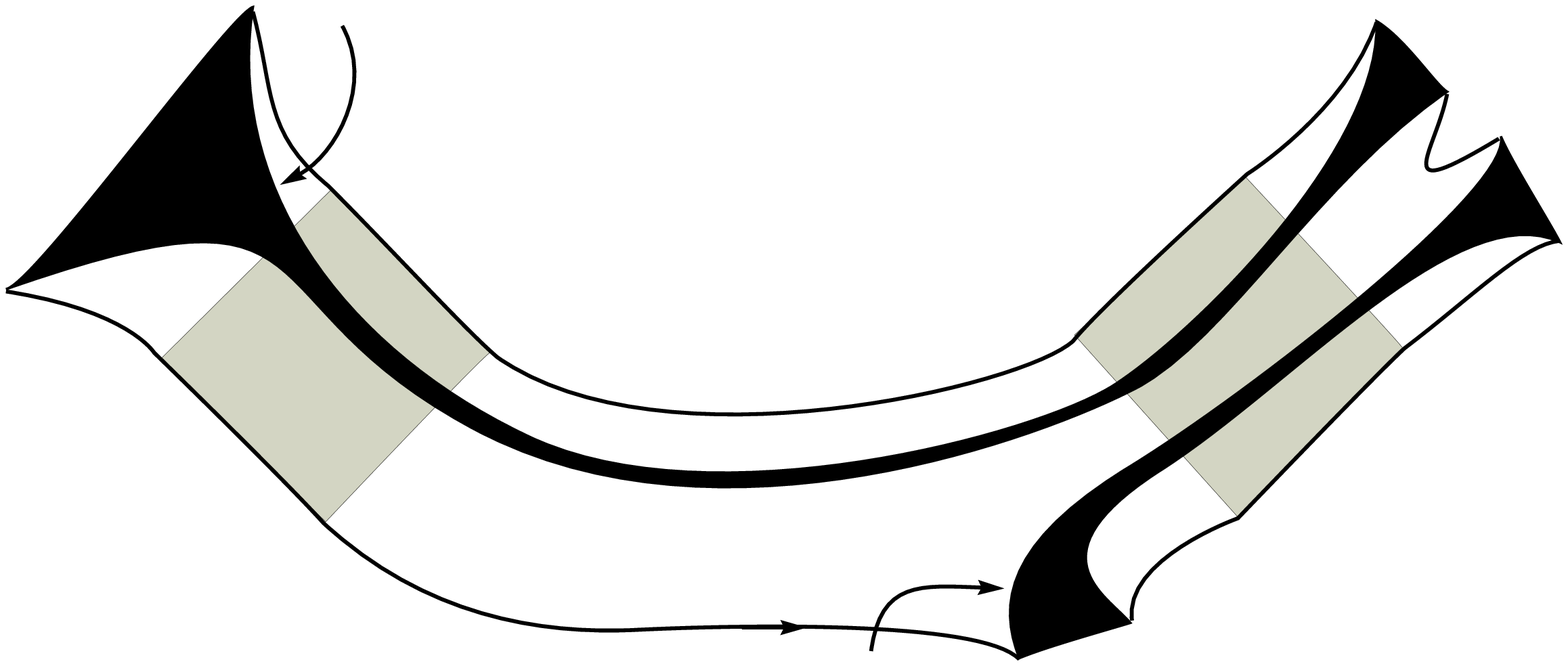}
\put (-170,125){\small$K$}
\put (-272,135){\small$B$}
\put (-305,158){\small$h (C)$}
\put (-250,188){\small$h (D_{2})$}
\put (-238,102){\small$A$}
\put (-100,102){\small$D$}
\put (-90.5,130){\small$D_{2}$}
\put (-160,70){\small$h (D_{1})$}
\put (-120,70){\small$h (B)$}
\put (-70,144){\small$C$}
\put (-35,158){\small$h (A)$}
\put (-54,181){\small$h (D)$}
\put (-250,126){\small$D_{1}$}

\end{picture}
\caption{The core of the leaf $L$ in Example~\ref{extwo}}\label{(2,2,2)core1}
\end{center}
\end{figure}

\begin{example}\label{extwo}

Carrying out the disk decomposition of $M$ (as in
Figure~\ref{(2,2,2)}) using the disks $+D_1,-D_2$ to obtain a
foliation $\FF$ of $M$, one can piece together the core of a leaf $L$
of $\FF$ as displayed in Figure~\ref{(2,2,2)core1} (the entire leaf
$L$ is displayed in Figure~\ref{leafL}(b) with the core shaded). Then
$$ A = \begin{bmatrix} 0&1\\ 1&1
\end{bmatrix}.$$
That is, in this case, $1$ can follow $2$ and $2$ can follow $1$ or
$2$, but $1$ cannot follow $1$. The loops of minimal period are
\begin{align*}
\Gamma _{(\dots 1212 \dots)} &= \alpha_1 -\alpha_2\\ \Gamma _{(\dots
2222 \dots)} &= -\alpha_2
\end{align*}
 The foliation cone to which $\FF$ belongs is defined by the
inequalities \begin{align*}
\{\alpha _{1}\ge \alpha _{2}\} &= \text{the half plane below the line
$\alpha _{1}=\alpha _{2}$}\\ \{\alpha _{2}\le0\} &= \text{the lower
half plane}. 
\end{align*}
  A base for this cone is the segment $[{\mathbf{e}_0},{\mathbf{e}_1}]$ (see
Figure~\ref{(2,2,2)cones}).
\end{example}

\begin{example}\label{exthree}

If one decomposes the sutured manifold $M$ of Figure~\ref{(2,2,2)}
using the disk $-D_0$, one obtains a product sutured torus and thus a
depth one foliation $\FF$ of $M$. A leaf $L$ of this foliation and the
disk $D_0$ is given in Figure~\ref{onedisk}. If $t$ is translation to
the left by one unit and $d$ is a Dehn twist in the dotted curve, then
$h=d^2\circ t$ is the pseudo-Anosov endperiodic monodromy of the
foliation $\FF$. The foliation $\FF$, the leaf $L$ and the endperiodic
monodromy are the same as in Example~\ref{exone} (see
Figure~\ref{leafL}(a)). The disk $D_0$ does contain $Z=\ZZ\cap L$ but
is not a Markov partition for $Z$. However the two components of
$h(D_0)\cap D_0$ will be a Markov partition for $Z$ and can be used to
compute the foliation cone with base $[{\mathbf{e}_1},{\mathbf{e}_2}]$ much
as above.

\begin{figure}[ht!]
\begin{center}
\begin{picture}(300,80)(10,100)

\epsfxsize=320pt
\epsffile{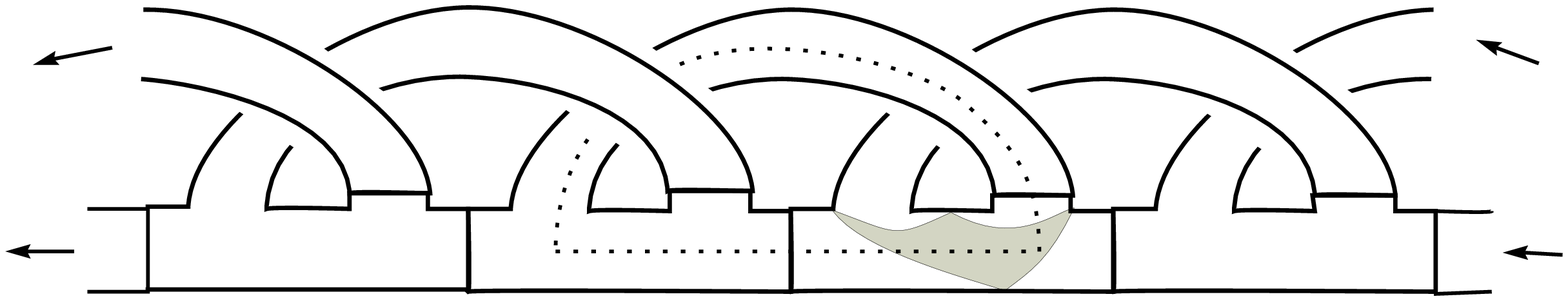}

\end{picture}
\caption{The leaf $L$ of Example~\ref{exthree} showing the disk $D_0$
(shaded)}\label{onedisk} 
\end{center}
\end{figure}

\end{example}

\begin{rem}
There are 13 knots of ten or fewer crossings that have unique Seifert
surface $S$ and whose complements, cut apart along the Seifert
surface, reduce via disks meeting the sutures twice to the sutured
manifold $M$ in Figure~\ref{(2,2,2)}(b).  They are: $$8_{15}, 9_{25},
9_{39}, 9_{41}, 9_{49}, 10_{58}, 10_{135}, 10_{144}, 10_{163},
10_{165}, 10_{49}, 10_{66}, 10_{80}$$ (with notation as in
Rolfsen~\cite{rolf}).  If $\kappa$ is one of the first ten of these
knots and $M$ is the sutured manifold of Figure~\ref{(2,2,2)}(b), then
there is a disk decomposition of $M_S(\kappa)$ to $M$ using two disks,
each of which meet the sutures twice.  By Proposition~\ref{twice}, the
foliation cones for $M_S(\kappa)$ are obtained from the foliation
cones in Figure~\ref{(2,2,2)cones} by crossing with $\R^2$. If
$\kappa$ is one of the last three of these knots, then there is a disk
decomposition of $M_S(\kappa)$ to $M$ using four disks, each of which
meet the sutures twice, and the foliation cones for $M_S(\kappa)$ are
obtained from the foliation cones in Figure~\ref{(2,2,2)cones} by
crossing with $\R^4$.  Finally, an additional knot $10_{53}$ has two
disjoint, nonisotopic Seifert surfaces $S_{1}\text{ and }S_{2}$. If we
cut the knot complement apart along $S_{1}\cup S_{2}$, we get two
disjoint sutured manifolds $X\text{ and }Y$.  The component $X$ is the
sutured manifold $X (2;1,2)$ defined in~\cite[page 385]{cc:surg} and
$Y=M_S(8_{15})$.  By results in~\cite{cc:surg}, the foliation cones in
$H^{1} (X)=\R^{4}$ consist of two half spaces, while those in $H^{1}
(Y)=\R^{4}$ were determined above. Using Mayer--Vietoris as
in~\cite[section 4]{cc:surg}, one shows that $$ \R^{4}=H^{1} (M_{S_{i}}
(10_{53}))\hra H^{1} (X)\oplus H^{1} (Y),\quad i=1,2, $$ is a diagonal
inclusion sending each foliated class to the direct sum of foliated
classes.  One then assembles the cone structure from that for $X\text{
and }Y$.
\end{rem}

\begin{figure}[ht!]
\begin{center}
\begin{picture}(300,150)(10,220)

\epsfxsize=320pt 
\epsffile{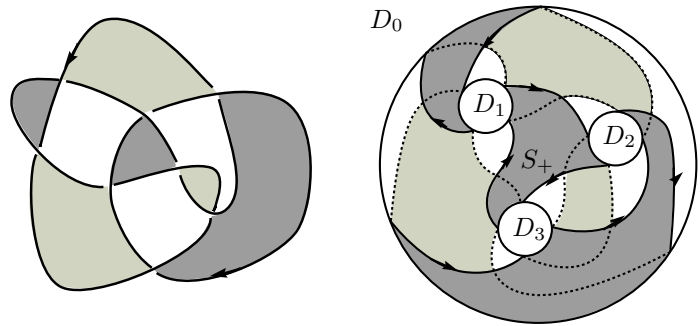} 
\put (-155,338){\small$D_{0}$} 
\put(-116,306){\small$D_{1}$} 
\put (-67,294){\small$D_{2}$} 
\put(-101,259){\small$D_{3}$} 
\put (-98,285){\small$S_{+}$}

\end{picture}
\caption{The $2$--component link of Example~\ref{exfour} and the
sutured manifold obtained from this link}\label{K2} 
\end{center}
\end{figure}

\begin{example}\label{exfour}

Consider the $2$--component link with Seifert surface in
Figure~\ref{K2} and the sutured manifold $M$ obtained from the link
complement by cutting apart along the Seifert surface.  We let
$D_1,D_2,D_3$ be the three disks indicated in the figure, $D_0$ the
outside disk. Then $H_2(M,\bd M) =\R^3$. The Thurston ball for $M$ and
bases for the foliation cones are given in Figure~\ref{K2cones}. One
computes the Thurston ball by noting that the given vertices all have
norm one and, in this case, the cones over the top dimensional faces
of the Thurston ball are all of the form $-C_i\cap C_j$ where
$C_i\text{ and }C_j$ are foliation cones.

\begin{figure}[ht!]
\begin{center}
\begin{picture}(300,145)(0,200)

\epsfxsize=300pt
\epsffile{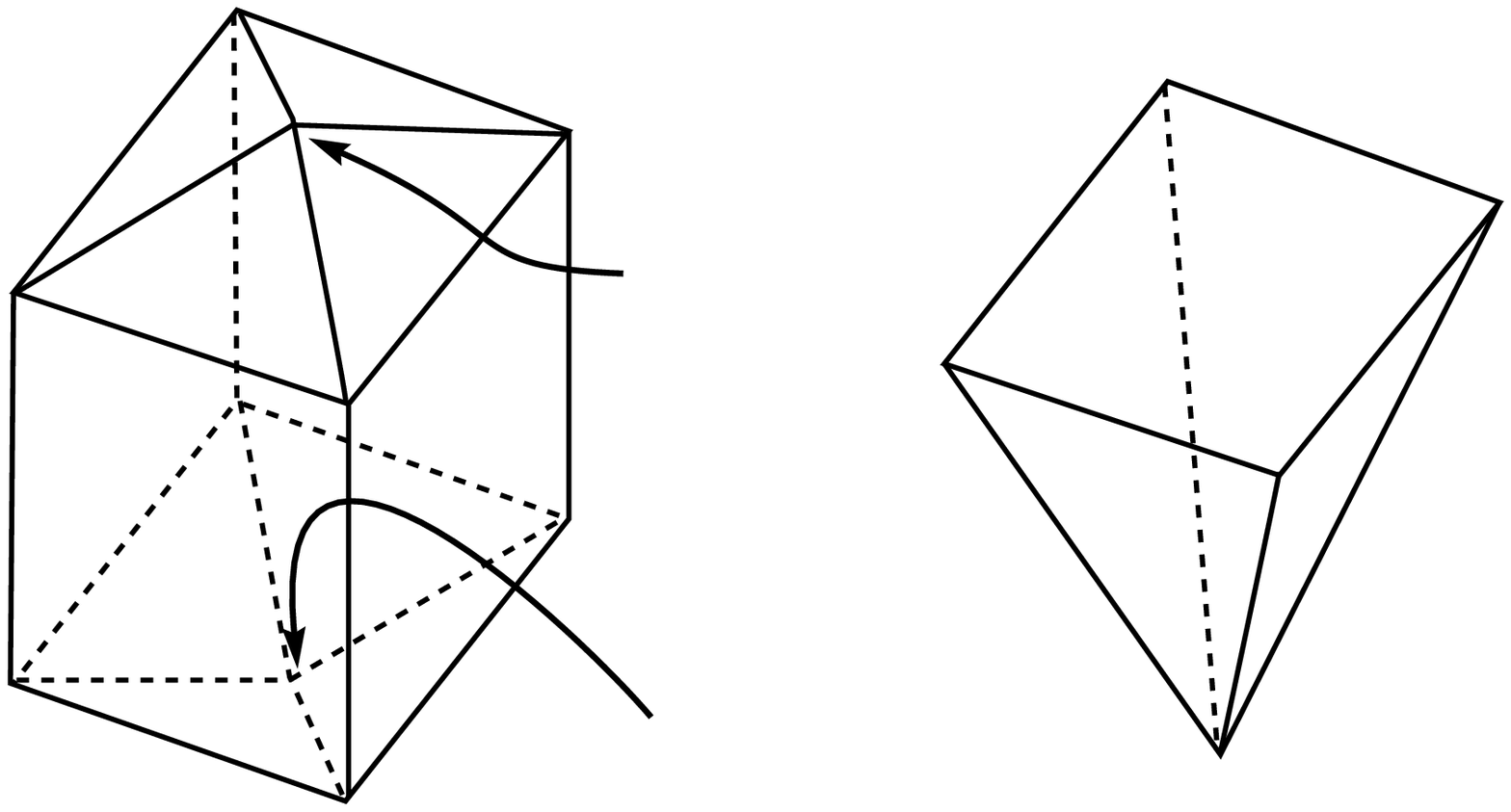}
\put (-210,205){\small$\mathbf{e}_{1}$}
\put (-270,222){\small$-\mathbf{e}_{2}$}
\put (-247,198){\small (a)}
\put (-275,287){\small$\mathbf{e}_{3}$}
\put (-127,280){\small$\mathbf{e}_{3}$}
\put (-242,340){\small$-\mathbf{e}_{1}$}
\put (-95,328){\small$-\mathbf{e}_{1}$}
\put (-240,275){\small$\mathbf{e}_{0}$}
\put (-208,273){\small$-\mathbf{e}_{0}$}
\put (-75,270){\small$-\mathbf{e}_{0}$}
\put (-173,316){\small$\mathbf{e}_{2}$}
\put (-25,306){\small$\mathbf{e}_{2}$}
\put (-173,250){\small$-\mathbf{e}_{3}$}
\put (-179,218){\small$-(\mathbf{e}_{2}+\mathbf{e}_{3})$}
\put (-95,212){\small$-(\mathbf{e}_{2}+\mathbf{e}_{3})$}
\put (-94,198){\small (b)}
\put (-165,295){\small$(\mathbf{e}_{2}+\mathbf{e}_{3})$}

\end{picture}
\caption{(a)  The Thurston ball and (b) the foliation cones for
Example~\ref{exfour}}\label{K2cones}
\end{center}
\end{figure}

\begin{figure}[ht!]
\begin{center}
\begin{picture}(300,150)(-10,45)

\epsfxsize=290pt
\epsffile{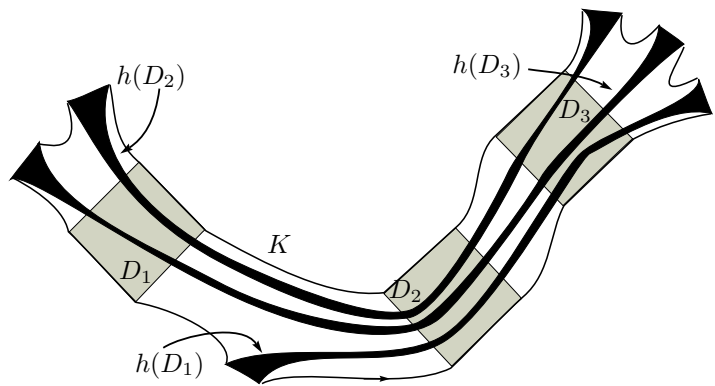}
\put (-180,95){\small$K$}
\put (-134,78){\small$D_{2}$}
\put (-236,85){\small$D_{1}$}
\put (-237,159){\small$h(D_{2})$}
\put (-110,163){\small$h(D_{3})$}
\put (-70,146){\small$D_{3}$}
\put (-230,50){\small$h(D_{1})$}

\end{picture}
\caption{The core $K$ of the leaf  in Example~\ref{exfour}}\label{K2core}
\end{center}
\end{figure}

We compute the foliation cone with square base in
Figure~\ref{K2cones}(b). If one carries out the disk decomposition of $M$
using the disks $-D_1,+D_2,+D_3$ to obtain a foliation $\FF$ of $M$,
one can piece together the core of a depth one leaf $L$ of $\FF$ as
displayed in Figure~\ref{K2core}. Then $$ A= \begin{bmatrix} 0&1&1\\
1&1&1\\ 1&1&1
\end{bmatrix}.
$$ That is, in this case, $1$ can follow $2$ and $3$, while $2$ and
$3$ can each follow $1$, $2$ and $3$. The minimal period loops are
\begin{align*}
\Gamma _{(\dots 1212
\dots)} &= \alpha_2 - \alpha_1\\ \Gamma _{(\dots 2222 \dots)} &=
\alpha_2\\  \Gamma _{(\dots 1313 \dots)}&= \alpha_3 - 
\alpha_1\\  \Gamma _{(\dots
3333 \dots)}&=\alpha_3 .
\end{align*}
 The foliation cone with square base   is
defined by the inequalities \begin{align*}
\alpha _{2}&\ge \alpha _{1}\\ \alpha _{2}&\ge0\\\alpha _{3}&\ge \alpha
_{1}\\ \alpha _{3}&\ge0.
\end{align*}

If $\kappa = 9_{38}$ or $10_{97}$ and $M$ is the sutured manifold in
Figure~\ref{K2}, then one can do a disk decomposition of $M_S(\kappa)$
to $M$ using a disk that meets the sutures twice. By
Proposition~\ref{twice}, the foliation cones for $M_S(\kappa)$ are
obtained from the foliation cones in Figure~\ref{K2cones}(b) by
crossing with $\R$.
\end{example}

\begin{example}\label{hex}
If $\kappa = 10_{55}$ or $10_{63}$ and $M$ is the sutured manifold
consisting of the complement of the pretzel link $(2,4,2)$ cut apart
along the Seifert surface (as in Figure~\ref{224}), then one can do a
disk decomposition of $M_S(\kappa)$ using two disks, each disk meeting
the sutures twice, obtaining $M$. By Proposition~\ref{twice}, the
foliation cones for $M_S(\kappa)$ are obtained from the foliation
cones in Figure~\ref{224} by crossing with $\R^2$.  Figure~\ref{224}
also gives the Thurston ball (dashed) for $M$. As before if one
carries out a disk decomposition of $M$ with the disks $+D_0,-D_1$, it
is easy to trace out a suitable core of the leaf
$L$. Figure~\ref{224core} gives this core as well as the rectangles
$R_{i}$ used in the Markov partition. The disk $D_{1}$ is represented
as hexagonal since it meets the sutures $6$ times.  The rectangular
disk $D_{0}=R_{1}$ meets the sutures $4$ times.  Notice that $D_{1}$
is divided into three rectangles $R_{2},R_{3},\text{ and }R_{4}$ of
the partition.  These look like pentagons, but the barycenter of
$D_{1}$ is not to be viewed as a vertex of these rectangles. The
incidence matrix for this scheme is $$ A =
\begin{bmatrix} 0&0&1&0\\ 1&0&1&0\\ 1&0&0&1\\ 0&1&0&0
\end{bmatrix}.
$$ 
The  minimal period loops are 
\begin{align*}
\Gamma _{(\dots 234234 \dots)} &= \alpha_2 - \alpha_1\\ \Gamma
_{(\dots 131313 \dots)} &= (-\alpha_2) + (\alpha_2 - \alpha_1) =
-\alpha_1.
\end{align*}
 The foliation cone is defined by the inequalities  \begin{align*}
\alpha _{2}&\ge \alpha _{1}\\ \alpha _{1}&\le0
\end{align*}
 and has base $[{\mathbf{e}_0},{\mathbf{e}_2}]$.

\begin{figure}[ht!]
\begin{center}
\begin{picture}(300,150)(0,60)

\epsfxsize=300pt
\epsffile{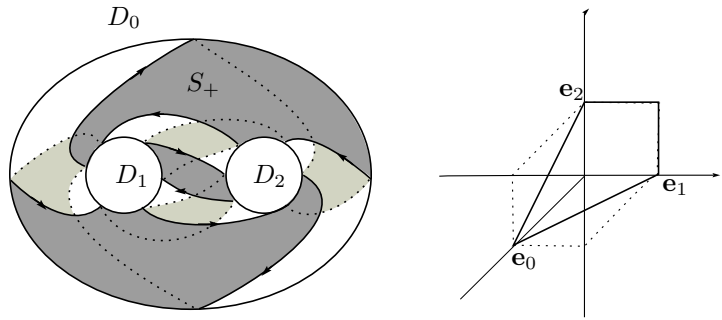}
\put (-250,180){\small$D_{0}$}
\put (-247,120){\small$D_{1}$}
\put (-195,120){\small$D_{2}$}
\put (-220,155){\small$S_{+}$}
\put (-40,117){\small$\mathbf{e}_{1}$}
\put (-97,89){\small$\mathbf{e}_{0}$}
\put (-79,153){\small$\mathbf{e}_{2}$}

\end{picture}
\caption{The sutured manifold, foliation cones, and Thurston ball for
$(2,4,2)$}\label{224} 
\end{center}
\end{figure}

\begin{figure}[ht!]
\vglue-0.3in
\begin{center}
\begin{picture}(300,130)(0,70)

\epsfxsize=300pt
\epsffile{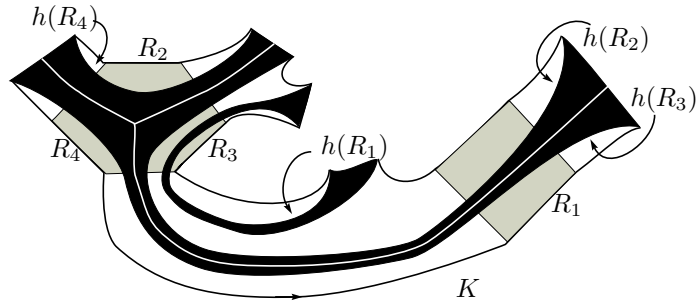}
\put (-100,75){\small$K$}
\put (-64,107){\small$R_{1}$}
\put (-34,147){\small$h (R_{3})$}
\put (-52,170){\small$h (R_{2})$}
\put (-194,128){\small$R_{3}$}
\put (-151,129){\small$h (R_{1})$}
\put (-253.5,128){\small$R_{4}$}
\put (-220.5,166){\small$R_{2}$}
\put (-260.5,178){\small$h(R_{4})$}

\end{picture}
\caption{The core $K$ of the leaf in Example~\ref{hex} }\label{224core}
\end{center}
\end{figure}

\end{example}

\begin{figure}[ht!]
\vglue -0.3in
\begin{center}
\begin{picture}(300,150)(0,60)

\epsfxsize=300pt
\epsffile{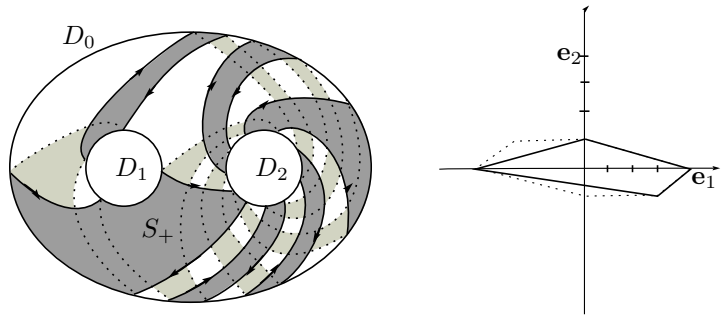}
\put (-247,120){\small$D_{1}$}
\put (-237,97){\small$S_{+}$}
\put (-267,170){\small$D_{0}$}
\put (-194,120){\small$D_{2}$}
\put (-80,164){\small$\mathbf{e}_{2}$}
\put (-29,117){\small$\mathbf{e}_{1}$}

\end{picture}
\caption{Example in which $\{D_1,D_2\}$ does not work as Markov
partition}\label{dodec} 
\end{center}
\end{figure}

\begin{rem}
In Example~\ref{hex} one is tempted to take
$\{R_{1}'=D_0,R_{2}'=D_1\}$ as ``Markov partition''. The hexagon
$h(R_{2}')$ stretches and passes properly through both $R_{2}'$ and
$R_{1}'$, while $h(R_{1}')$ meets only $R_{2}'$. In fact, the
resulting $2\times 2$ incidence matrix does determine the correct
cone.  In general, however, one must stick to the standard definition
of ``Markov partition'', requiring that the elements be rectangles
with one pair of opposite edges in the stable set and one pair of
opposite edges in the unstable set. This is illustrated by the sutured
manifold in Figure~\ref{dodec}.  We omit the details, but it happens
in this case that, carrying out the disk decomposition by
$-D_{1}\text{ and} -\!D_{2}$ and using $\{D_1,D_2\}$ as ``Markov
partition'', one does not get the entire cone with base
$[-{\mathbf{e}_1},(3{\mathbf{e}_1}-{\mathbf{e}_2})/4]$.  In this case,
$D_{2}$ is a decagon.  Figure~\ref{dodec} also gives the bases for the
foliation cones and the Thurston ball (dashed) for this example.

\end{rem}

\begin{rem}

Of the 249 knots of ten or fewer crossings, 117 are fibered and 111
belong to the class $\BB$ of knots described in~\cite{cc:surg}. The
preceding examples indicate how to compute the foliation cones of 18
of the remaining 21 knots. The three remaining knots,
$9_{35},10_{101},10_{120}$, of $\le10$ crossings each have a unique
Seifert surface. Let $M$ be the knot complement cut apart along this
Seifert surface. The knot $9_{35}$ is the pretzel knot $(3,3,3)$, $M$
is completely reduced and $H_2(M,\bd M) = \R^2$. For the knot
$10_{101}$, decomposition by a disk meeting the sutures twice replaces
$M$ with a completely reduced, sutured manifold $M'$ with $H_2(M',\bd
M') =
\R^3$. For the knot $10_{120}$, $M$ is completely reduced and $H_2(M,\bd M) =
\R^4$. In all cases the foliation cones are
easily computed as above.

\end{rem}

\Addresses\recd

\newpage

\erratum1
\pagenumbers{571}{576}
\papernumber{27}
\received{18 August 2000}
\published{30 August 2000}
\let\reviseddate\relax

\newcommand{\OO}{\mathcal{O}}
\newtheorem{claim}{Claim}

\title{Foliation Cones;\\A Correction}\shorttitle{Erratum: Foliation cones}

\begin{abstract}
The proof of~\cite[Lemma~3.6]{cc:cone} was incorrect.
Happily, a correct proof can be given.
\end{abstract}
\asciiabstract{}

\primaryclass{57R30}
\secondaryclass{57M25, 58F17} 


\maketitle

The goal of~\cite[Section~3]{cc:cone} was to prove that Theorem~1.1
(the main result of the paper) holds for arbitrary sutured, depth--one
foliated manifolds $M$ provided that it holds for completely reduced
ones such that $\tb M$ has no toral or annular components. This was
essential for the application of the Handel--Miller theory in
Section~5.  The proof was by induction on the cardinality $r$ of a
maximal reducing family $\{T_{1},\dots,T_{r}\}$ for $(M,\tb
M)$~\cite[Definition~3.2]{cc:cone}. This number will be called the
``reducing rank''.

\begin{lemma}
At each step of the induction, no generality is lost by assuming that
each component of $\tb M$  is neither a torus nor an annulus.
\end{lemma}
 
This is essentially Lemma~3.6 of ~\cite{cc:cone} and was proven by
pointing out that, in either case, a small perturbation makes $\FF$
transverse to the toral or annular component $N$ of $\tb M$ and does
not increase the reducing rank.  In the case that $N$ is a torus, this
simply makes $N$ a component of $\trb M$ without changing the foliated
cohomology class.  There is no problem here as the resulting foliated
manifold remains tautly foliated and sutured. In the case of an
annulus, however, the resulting foliation is not taut and $(M,\gamma
)= (M,\trb M)$ does not satisfy the definition of a sutured manifold.
The correct strategy for avoiding annular components of $\tb M$
depends on a strengthening of Lemma~3.5 of~\cite{cc:cone}.

Recall that, in the statement of Lemma~3.5, a connected manifold $M'$
is obtained from $M$ by excising the interior of a normal neighborhood
$$N (T)=T\times [-1,1]$$ of a reducing surface $T\ss M$, $T$ being a
properly imbedded torus or annulus.  When $T$ is an annulus, it is
required that one component of $\bd T$ lie on an inwardly oriented
component of $\tb M$, the other on an outwardly oriented one.  The
strengthened version of Lemma~3.5 follows.

\begin{lemma}
If $M'$ is connected, the conclusion of {\em Theorem~1.1} holds for
$M'$ if and only if it holds for $M$.
\end{lemma}

The ``only if'' part of this lemma is the content
of~\cite[Lemma~3.5]{cc:cone}. The proof was omitted, but here we will
give the complete proof of Lemma~2.  First we use it as follows.

\begin{proof}[Proof of Lemma 1]
The only problem is in the case that a component $N$ of $\tb M$ is an
annulus. In this case, $N$ is flanked by two annular sutures
(components of $\trb M$) on which $\FF$ induces the foliation by
spirals.  Gluing these together so as to match the foliations produces
a tautly foliated, sutured manifold $(M^{*},\FF^{*})$, having a toral
component of $\tb M^{*}$ and one new properly imbedded annular
reducing surface $T_{r+1}$. The maximal reducing family is now
$\{T_{1},\dots,T_{r},T_{r+1}\}$.  Furthermore, $(M,\FF)$ is obtained
from $(M^{*},\FF^{*})$ by removing an open normal neighborhood of
$T_{r+1}$. By Lemma~2, Theorem~1.1 will hold for $M$ if it holds for
$M^{*}$.  If we perturb $\FF^{*}$ to be transverse to the new toral
component of $\tb M^{*}$, the annulus $T_{r+1}$ is no longer a
reducing surface. Any of the $T_j$, $1\le j\le r$, with one boundary 
on $N$ will no longer be a reducing surface, and no additional 
reducing surfaces will have been introduced. Hence, the maximal reducing 
family  in the new foliated manifold is a (possibly proper) subset 
of $\{T_{1},\dots,T_{r}\}$, but the tangential boundary has
one less annular component.  Finitely many repetitions of this
procedure removes all annular components but does not affect the
induction on the reducing rank.
\end{proof}

  The proof of Lemma~2 will use the
Mayer--Vietoris sequence \begin{equation}\tag{$*$}
\hr \xrightarrow{i} H^{1}(M';\R)\oplus H^{1}(N(T);\R) \xrightarrow{j}
H^{1}(N_{-};\R)\oplus H^{1}(N_{+};\R), \end{equation} where
\begin{align*} N_{-} 
&= T\times [-1,-1/2]\\ N_{+} &= T\times [1/2,1]. \end{align*} Denote
by $\OO (M)\text{ and }\OO (M')$ the open sets of foliated classes in
$H^{1} (M;\R)$ and $H^{1} (M';\R)$ respectively. Hereafter, we will
omit the coefficient field $\R$ from the notation for homology and
cohomology.

\begin{claim}
The Poincar\'e dual $[\alpha _{T}]$ of $[T]\in H^{2} (M,\bd M)$ spans
$\ker i=\ker \varphi $, where $\varphi:H^{1} (M)\to H^{1} (M') $ is
the restriction map.
\end{claim}

\begin{proof} 
Examining the degree 0 terms of the Mayer--Vietoris sequence reveals
that $i$ has 1--dimensional kernel.  The Poincar\'e dual of $[T]$ is
represented by a closed 1--form $\alpha_{T} $, chosen to have compact
support in $T\times (-1/2,1/2)$. Since $T$ does not separate $M$,
$0\ne [\alpha_{T} ]\in\hr$.  Since $T$ does separate $N(T)$,
$(0,0)=([\alpha_{T} |M'],[\alpha_{T} |N(T)])$ and the assertion for $\ker i$
follows. If $$p:H^{1}(M';)\oplus H^{1}(N(T);)\ra H^{1}(M';)$$ is
projection onto the first summand, we write $\varphi =p\o i$ and prove
that $p$ is one--to--one on $\im i$. By exactness of $(*)$, $\im
i=\ker j$, so an element of $\im i$ annihilated by $p$ must be of the
form $(0,[\eta ])$, where $$ 0 = j (0,[\eta ]) = (-[\eta
|N_{+}],-[\eta |N_{-}]).  $$ Since the inclusions $N_{\pm}\hra N (T)$
are homotopy equivalences, $[\eta ]=0$.
\end{proof}

\begin{claim}
For $\varphi :H^{1} (M)\to H^{1} (M')$ as above, \begin{gather*}
\varphi (\OO (M)) = \OO (M')\\ \varphi ^{-1} (\OO (M')) = \OO (M).
\end{gather*}
\end{claim}

\begin{proof}
If $[\omega ]\in\OO (M)$, we assume that $\omega $ is a foliated form
and let $\FF$ be the corresponding foliation.  By a theorem of
Roussarie and Thurston~\cite{rous:plongenew,th:normnew,cc:rt}, an isotopy
moves $T\times \{\pm1\}$ to a position everywhere transverse to $\FF$
(fixing $\bd T\times \{\pm1\}$ in the case that $T$ is an annulus).
Thus, $[\omega |M']=\varphi [\omega ]$ is a foliated class.

Now let $[\omega ']\in\OO (M')$, where $\omega '$ is a foliated form.
Let $\FF'$ be the foliation defined by $\omega '$.  Since
$T_{+}=T\times
\{1\}$ and $T_{-}= T\times \{-1\}$ are homologous in $M'$, we see that
the forms $\omega '|T_{\pm}$ induce cohomologous forms on $T$.  If $T$
is an annulus, the foliations $\FF|T_{+}$ and $\FF|T_{-}$ are either
both product foliations or both foliations by spirals.  In either
case, $\FF'$ extends across $N\times [-1,1]$ to provide a taut
foliation of $M$ with foliated class $[\omega ]$ such that $\varphi
[\omega ]=[\omega ']$.  In the case that $T$ is a torus, the theorem
of Laudenbach and Blank~\cite{LBnew} implies that, after an isotopy in $M'$
supported near $T_{\pm}$, $\FF'$ again extends across $N\times [-1,1]$
and provides the foliated class $[\omega ]$ such that $\varphi
[\omega ]=[\omega ']$.  Together with the previous paragraph, this
proves that \begin{gather*}
\varphi (\OO (M)) = \OO (M')\\ \OO (M) \sseq \varphi ^{-1} (\OO (M')).
\end{gather*}

Finally we must show that, if $\varphi [\omega ]=[\omega ']\in\OO
(M')$, then $[\omega ]\in\OO (M)$.  Write $i[\omega ]= ([\omega
'],[\eta ])$ and choose the representative form $\omega '$ to be
foliated. After an isotopy, we can also assume that $T\times \{1\}$ is
transverse to $\omega '$.  Exactness of the sequence~$(*)$ implies
that $\omega '\text{ and }\eta $ restrict to cohomologous forms on
$N_{\pm}$. Replacing $\eta $ with a cohomologous form $\eta +df$, we
can assume without loss of generality that $\eta |N_{+}=\omega'
|N_{+}$. By an isotopy within $N(T)$, we compress $N(T)$ into a
neighborhood of $N_{+}$ where $\eta $ is nonsingular, keeping $N_{+}$
itself pointwise fixed and carrying $T\times \{-1\}$ to a position
transverse to $\eta $.  Reversing this isotopy, we see that no
generality is lost in assuming that $\eta $ is a foliated form on
$N(T)$.  As above, there is an isotopy of $\eta $, supported in a
neighborhood of $N_{-}$ in $N(T)$, to a form $\eta '$ agreeing with
$\omega '$ on $N_{\pm}$.  Then $\omega '\text{ and }\eta '$ assemble
to a foliated form $\wt{\omega }$ on $M$ and $$i([\wt{\omega}
])=([\omega '],[\eta '])=i([\omega ]).$$ By Claim~1, $[\omega
]=[\wt{\omega }]+c[\alpha_{T} ],$ where $c\in\R$ and the class
$[\alpha_{T} ]$ is Poincar\'e dual to $[T]=[T\times
\{-1\}]$.  Since $\wt{\omega }$ is transverse to $T\times \{-1\}$, we
can choose $\alpha_{T} $ to be compactly supported near $T\times
\{-1\}$ and to vanish identically on a vector field $v$ (on $M\sm\tb
M$) such that $\wt{\omega }(v)>0$ everywhere.  That is, the closed
form $\wh{\omega }=\wt{\omega }+c
\alpha_{T} $ is nonsingular.  Also, $\alpha_{T} $ is bounded, implying that
$\wh{\omega }$ also blows up nicely at $\tb M$ and $[\omega
]=[\wh{\omega }]$ is a foliated class.
\end{proof}

Let $\ell_{T}\ss H^{1} (M)$ be the one--dimensional subspace spanned
by  the class $[\alpha _{T}]$.

\begin{claim}
If $U$ is a connected component of $\OO (M)$, then $\ell_{T}\ss\ol{U}$.

\end{claim}

\begin{proof}
Indeed, let $\omega $ be a foliated form representing an element of
$U$. As in the previous proof, $\alpha _{T}$ can be chosen to vanish
identically on a vector field $v$ such that $\omega (v)>0$ everywhere.
Thus, for each $c[\alpha _{t}]\in\ell_{T}$, the line segment $$
t[\omega ] + (1-t)c[\alpha _{T}], \quad0\le t\le1, $$ connects
$c[\alpha _{T}]$ to $[\omega ]$ and lies in $U$ for $t>0$.
\end{proof}

\begin{proof}[Proof of Lemma 2]
First we eliminate trivial cases.  By Claim~2, $0\in\OO (M)$ if and
only if $0\in\OO (M')$.  By~\cite[Proposition~3.7]{cc:cone}, $0$ is a
foliated class if and only if the manifold is a product $S\times I$ of
a compact surface $S$ and a compact interval $I$.  In turn, this is
the case if and only if the entire cohomology space is the unique
foliation cone and satisfies Theorem~1.1 of~\cite{cc:cone} trivially.
Thus, we assume that neither $M$ nor $M'$ is a product. Claim~2 also
allows us to assume that neither $\OO (M)$ nor $\OO (M')$ are empty.

Since $\OO (M')$ is open in the vector space $H^{1} (M')$, Claim~2
implies that the linear map $$\varphi :H^{1} (M)\to H^{1} (M')$$ is
surjective. If Theorem~1.1 holds on $M'$, we can use this surjection
to pull the cone structure back to $H^{1} (M)$ and use Claim~2 to
verify that Theorem~1.1 holds on $M$. For the converse, suppose the
theorem holds on $M$.  If $\CC\ss H^{1} (M)$ is a foliation cone, the
fact that it is neither empty nor the entire vector space implies that
it is defined by a finite set of nontrivial linear inequalities
$\theta _{i}\ge0$, $1\le i\le q$.  By Claim~3, $\ell_{T}\ss\CC$, hence
$\theta _{i}|\ell_{T}\ge 0$ and this implies that $\theta
_{i}|\ell_{T}\equiv0$.  By Claim~1, the linear functionals $\theta
_{i}$ pass to nontrivial linear functionals $\wt{\theta }_{i}$ on
$H^{1} (M')$.  The convex, polyhedral cone $\CC'$ defined by the
linear inequalities $\wt{\theta }_{i}\ge 0$ is precisely the image of
$\CC$ under $\varphi $ and has $\CC$ as its entire pre--image.  By
Claim~2, Theorem~1.1 follows easily for $M'$.
\end{proof}

\Addresses\recd


\begin{thebibliography}



\bibitem{bca}
{\bf S\,A Bleiler}, {\bf A\,J Casson}, {\em Automorphisms of surfaces
after {N}ielsen 
  and {T}hurston}, Cambridge Univ. Press (1988)

\bibitem{cc:LB} {\bf J Cantwell}, {\bf L Conlon}, {\em Isotopies of
foliated $3$--manifolds without
  holonomy},  Adv. in Math. 144 (1999) 13--49

\bibitem{cc:prop2}
{\bf J Cantwell}, {\bf L Conlon}, {\em Leafwise hyperbolicity of
proper foliations}, Comment. Math. 
  Helv.  64 (1989) 329--337

\bibitem{cc:isotopy}
{\bf J Cantwell}, {\bf L Conlon}, {\em Isotopy of depth one
foliations}, {\it Proceedings of the 
  {I}nternational {S}ymposium and {W}orkshop on Geometric Study of Foliations,
  Tokyo, 1993}, World Scientific,  (1994) 153--173

\bibitem{cc:surg} {\bf J Cantwell}, {\bf L Conlon}, {\em Surgery and
foliations of knot complements}, Journal of Knot
  Theory and its Ramifications,  2 (1993) 369--397

\bibitem{cc:smth2}
{\bf J Cantwell}, {\bf L Conlon}, {\em Topological obstructions to
smoothing proper foliations}, 
  Contemporary Mathematics,  161 (1994) 1--20

\bibitem{dip}
{\bf P Dippolito}, {\em Codimension one foliations of closed
manifolds}, Ann. of Math.  107 (1978) 403--453

\bibitem{E-S} {\bf C\,J Earle}, {\bf A Schatz}, {\em Teichm\"uller
theory for surfaces with boundary}, J. Diff. Geo.  4 (1970) 169--185

\bibitem{fe:endp} {\bf S Fenley}, {\em Endperiodic surface
homeomorphisms and $3$--manifolds}, Math.  Z. 224 (1997) 1--24

\bibitem{F-O}
{\bf W Floyd}, {\bf H Oertel}, {\em Incompressible surfaces via
branched surfaces}, 
  Topology, 23 (1984) 117--125

\bibitem{fried} {\bf D Fried}, {\em Fibrations over ${S}^{1}$ with
pseudo-{A}nosov monodromy},
  Ast\'erisque,  66--67 (1991) 251--266

\bibitem{ga1} {\bf D Gabai}, {\em Foliations and the topology of
$3$--manifolds}, J. Diff. Geo.
   18 (1983) 445--503

\bibitem{ga0} {\bf D Gabai}, {\em Foliations and genera of links},
Topology, 23 (1984)
  381--394

\bibitem{gab:kneser} {\bf D Gabai}, {\em Essential laminations and
Kneser normal form}, preprint

\bibitem{hae}
{\bf A Haefliger}, {\em Vari\'et\'es feuilletees}, Ann. Scuola Norm.
Sup. Pisa,   16 (1962) 367--397

\bibitem{ham:disk-holes} {\bf M\,E Hamstrom}, {\em Some global
properties of the space of homeomorphisms on a
  disk with holes}, Duke Math. J.  29 (1962) 657--662

\bibitem{ham:torus} {\bf M\,E Hamstrom}, {\em The space of
homeomorphisms on a torus}, Ill. J. Math. 
  9 (1965) 59--65

\bibitem{ham:homeo}
{\bf M\,E Hamstrom}, {\em Homotopy groups of the space of homeomorphisms on a
  $2$--manifold}, Ill. J. Math.  10 (1966) 563--573

\bibitem{LB} {\bf F Laudenbach}, {\bf S Blank}, {\em Isotopie de
formes ferm\'ees en dimension
  trois}, Inv. Math.  54 (1979) 103--177

\bibitem{Oer} {\bf U Oertel}, {\em Incompressible branched surfaces},
Inv. Math.  76
  (1984) 185--410

\bibitem{Oer:hom}
{\bf U Oertel}, {\em Homology branched surfaces: {T}hurston's norm on ${H}_{2}
  ({M}^{3})$}, from: ``Low {D}imensional {T}opology and {K}leinian {G}roups'',
  (D\,B\,A
  Epstein, editor), Lond. Math. Soc. Lecture
  Notes,  112, Cambridge University Press (1985) 253--272

\bibitem{plante:meas} {\bf J\,F Plante}, {\em Foliations with measure
preserving holonomy}, Ann. of Math.
 102 (1975) 327--362

\bibitem{rolf}
{\bf D Rolfsen}, {\em Knots and Links}, Publish or Perish, Inc.
(1976) 

\bibitem{rous:plonge}
{\bf R Roussarie}, {\em Plongements dans les vari\'et\'es feuillet\'ees et
  classification de feuilletages sans holonomie}, I.H.E.S. Sci. Publ. Math.
 43 (1973) 101--142

\bibitem{sa:pseudo} {\bf R Sacksteder}, {\em Foliations and
pseudogroups}, Amer. J. Math.  87
  (1965) 79--102

\bibitem{scharle} {\bf M Scharlemann}, {\em Sutured manifolds and
generalized {T}hurston norms}, J.
  Diff. Geo.  29 (1989) 557--614

\bibitem{sch_cycles} {\bf S Schwartzmann}, {\em Asymptotic cycles},
Ann. of Math.  66 (1957)
  270--284

\bibitem{scott_homeo} {\bf G\,P Scott}, {\em The space of
homeomorphisms of a $2$--manifold}, Topology,
 9 (1970) 97--109

\bibitem{sull:cycles} {\bf D Sullivan}, {\em Cycles for the dynamical
study of foliated manifolds and
  complex manifolds}, Inv. Math.  36 (1976) 225--255

\bibitem{th:norm} {\bf W Thurston}, {\em A norm on the homology of
three--manifolds}, Mem. Amer. Math.
  Soc.  59 (1986) 99--130

\bibitem{wald} {\bf F Waldhausen}, {\em On irreducible $3$--manifolds
which are sufficiently
  large}, Ann. of Math.  87 (1968) 56--88


\end{thebibliography}

\begin{thebibliography}

\bibitem{cc:cone} {\bf J\,Cantwell}, {\bf L\,Conlon}, {\em Foliation
cones}, Geometry and Topology Monographs, Volume 2, Proceedings of the
Kirbyfest, (1999) 35--86


\bibitem{cc:rt} {\bf J\,Cantwell}, {\bf L\,Conlon}, {\em General
position in tautly foliated, sutured manifolds}, preprint

\bibitem{LBnew} {\bf F\,Laudenbach}, {\bf S\,Blank}, {\em Isotopie de
formes ferm\'ees en dimension trois}, Inv. Math. 54 (1979) 103--177

\bibitem{rous:plongenew} {\bf R\,Roussarie}, {\em Plongements dans les
vari\'et\'es feuillet\'ees et classification de feuilletages sans
holonomie}, I.H.E.S. Sci. Publ. Math. 43 (1973) 101--142

\bibitem{th:normnew} {\bf W\,Thurston}, {\em A norm on the homology of
three--manifolds}, Mem. Amer. Math. Soc. 59 (1986) 99--130

\end{thebibliography}
\end{document}